\documentclass[12pt]{article}
\date{}

\usepackage[a4paper, total={7.3in, 9in}]{geometry}

\usepackage{amsmath,mathrsfs}

\usepackage{shapepar,bm}
\usepackage[dvips,arrow,matrix,ps,color,line]{xy}
\usepackage{theorem,amsfonts,amssymb,amscd}
\usepackage{geometry}
\usepackage{makeidx}
\usepackage{eulervm}

\usepackage{graphicx}
\usepackage{graphics}

\usepackage{epstopdf}

\usepackage{amssymb,graphics,hyperref}

\hypersetup{
    colorlinks = true,
    linkcolor = blue,
    anchorcolor = red,
   citecolor = blue,
    filecolor = red,
    pagecolor = red,
    urlcolor = blue}

\newcommand{\mathsym}[1]{{}}

\usepackage[usenames]{color}
\definecolor{MyLightMagenta}{cmyk}{0.1,0.8,0,0.1}
\definecolor{MyDarkBlue}{rgb}{0.1,0,0.3}
\makeindex

\def\wb{[\bfb]}
\def\wbet{[\bm\beta]}

\def\NN{\mathbb N}

\def\bfz{{\mathbf z}}

\def\ovsig{\overline{\sigma}}

\def\bfb{{\mathbf b}}

\def\gln{{gl_n(\ZZ)}}
\def\glin{{gl_\infty(\ZZ)}}

\def\bfu{{\mathbf u}}
\def\bfv{{\mathbf v}}
\def\bfw{{\mathbf w}}

\def\Rcal{\mathcal R}

\def\Bcal{{\mathcal B}}

\def\ZZ{\mathbb Z}

\def\QQ{\mathbb Q}
\def\PP{\mathbb P}

\def\cocoa{{\hbox{\rm C\kern-.13em o\kern-.07em C\kern-.13em o\kern-.15em A}}}

\def\Dcal{\mathcal D}

\def\Acal{{\mathcal A}}

\def\bfu{{\bf u}}

\def\ovDcal{\overline{\Dcal}}

\def\End{\mathrm{End}}

\def\blamb{{\bm \lambda}}
\def\bmu{{\bm \mu}}

\def\Pcal{{\mathcal P}}

\def\w2M{\bigwedge^2M}

\def\wM{\bigwedge M}

\def\w{\wedge }
\def\bw{\bigwedge }

\protect
\protect
\protect\def\wMn{{\bigwedge M_n}}
\protect

\protect
\protect

\def\sra{\rightarrow}

\def\proof{\noindent{\bf Proof.}\,\,}
\def\qed{{\hfill\vrule height4pt width4pt depth0pt}\medskip}
\def\be{\begin{equation}}
\def\ee{\end{equation}}
\def\bclm{\begin{claim}}
\def\eclm{\end{claim}}
\def\beqn{\begin{eqnarray}}
\def\eeqn{\end{eqnarray}}
\def\beqn*{\begin{eqnarray*}}
\def\eeqn*{\end{eqnarray*}}

\theoremstyle{change}
\theorembodyfont{\rmfamily}

\newtheorem{claim}{}[section]

\def\rr{\large{r}}

\title{Schubert Derivations on the Infinite Wedge Power}

\author{Letterio Gatto \& Parham Salehyan  \thanks{Work sponsored by FAPESP,  Processo n. 2016/03161-3 and, partially, by Politecnico di Torino, Finanzia-\newline mento 
Diffuso della Ricerca; INDAM-GNSAGA e PRIN "Geometria delle Variet\`a Algebriche''.
\newline ${}$ \,\,\,\,\,\, 2010 MSC: 14M15, 15A75, 05E05, 17B69.  
\newline ${}$\,\,\,\,\,\,\, {\em Keywords and phrases:} Hasse-Schmidt Derivations on Exterior Algebras, Schubert Derivations on infini-
\newline te wedge powers; Fermionic space, vertex operators, bosonic vertex representation of Date-Jimbo-Kashiwara-\newline Miwa.}}

\begin{document}

\maketitle

\vspace{-14pt}
\hfill{\em  To Israel Vainsencher on the occasion of his seventieth birthday {}\hskip 30pt ${}$}
\abstract{\noindent 
The {\em Schubert derivation} is a distinguished Hasse-Schmidt derivation on the exterior algebra of a free abelian group, encoding the formalism of Schubert calculus for all Grassmannians at once. The purpose of this paper is to  extend the Schubert derivation to the infinite exterior power of a free $\ZZ$-module of infinite rank (fermionic Fock space). Classical vertex operators naturally arise from the {\em integration by parts formula}, that also recovers the generating function occurring in  the {\em bosonic vertex representation} of the Lie algebra $\glin$,  due to Date, Jimbo, Kashiwara and Miwa (DJKM). In the present framework,  the DJKM result will be interpreted as a limit case of the following general observation: the singular cohomology of the complex Grassmannian $G(r,n)$ is an irreducible representation of the Lie algebra of $n\times n$ square matrices.}

\smallskip
\noindent
%\keywords{\em Hasse-Schmidt Derivations on Exterior Algebras; Schubert Derivations; Infinite Grassmannians; Bosonic and Fermionic Fock Spaces; Date-Jimbo-Kashiwara-Miwa bosonic vertex representation of the alebra $gl_\infty$.}

%In this paper we interpret the celebrated bosonic vertex representation by Date-Jimbo-Kashiwara and Miwa  as the limit of the following general  fact: the singular cohomology of the Grassmannian of $r$-dimensional subspaces of the $n$ dimensional complex vector space is a module over the Lie algebra of the $n\times n$ matrices.   A Schubert derivation is a natural Hasse-Schmidt derivation on an exterior algebra encoding the formalism of Schubert Calculus not just on one Grassmannian, but on all at once. However the exterior algebra contains no exterior power of infinite degree (infinite wedge power). In this paper we extend the Schubert derivation on Fermionic Fock spaces (infinite wedge powers of a free abelian group of infinite rank),  and using the related formalism we provide a new deduction of the celebrated bosonic vertex reresentation of the Lie algebra of all the matrices of infinite size,  whose  all the entries are zero but finitely many, due to pioneristic work by Date, Jimbo, Kashiwara and Miwa.}

\tableofcontents

\section*{Introduction} 
\addcontentsline{toc}{section}{\,\,\,\,\,Introduction}

\claim{\bf The Goal.} Let $r,n\in\NN\cup\{\infty\}$  such that $r\leq n$. The main 
characters of this paper are i) the exterior algebra $\bw M_n$ of a free abelian 
group $M_n:=\bigoplus_{0\leq j<n}\ZZ b_j$ and ii) the cohomology ring $B_{r,n}$ of 
the Grassmann variety $G(r,n)$. By the latter we mean the following. If $0\leq r
\leq n<\infty$,  $B_{r,n}$ stands for  the usual singular cohomology ring  
$H^*(G(r,n),\ZZ)$ of the Grasmannian variety parametrizing $r$-dimensional 
subspaces of the complex $n$-dimensional vector space. If $n=\infty$, the ring
$B_r:=B_{r,\infty}$ will denote the cohomology of the ind-variety $G(r,\infty)$ 
(see e.g. \cite[p.~302]{bottu} or \cite{DimiPenkov,IgnaPenkov}), which is a 
polynomial ring $\ZZ[e_1,\ldots,e_r]$ in $r$-indeterminates. If $r=n=\infty$,  instead
$Gr(\infty):=G(\infty,\infty)$ is the ind-Grassmannian constructed e.g. in \cite[Section 3.3]{IgnaPenkov} or  the Sato's Universal Grassmann Manifold (UGM), as in e.g. \cite{satoUGM}. In this case $B:=B_{\infty,\infty}$ is the $\ZZ$-polynomial ring in infinitely many indeterminates. Let $\Bcal_{ij}\in \End_\ZZ(M_n)$ such that $\Bcal_{ij}(b_k)=b_i\delta_{jk}$,  and let
\be
\gln:=\bigoplus_{0\leq i,j<n} \ZZ\cdot \Bcal_{ij}\subseteq \End_\ZZ(M_n),
\ee
which is a Lie algebra with respect to the usual commutator. Clearly $\gln=\End_\ZZ(M_n)$ if $n<\infty$.
This paper is inspired  by the following simple observation,   for which  we have not been able to find an explicit reference in the literature:

\begin{quotation}
\begin{center}
\noindent
{\em The ring $B_{r,n}$ is a module over the Lie algebra $\gln$}.
\end{center}
\end{quotation}

\noindent
If $r=1$ the claim is obvious, because if $e$ denotes the hyperplane class of $\PP^{n-1}$, then  $B_{1,n}=\ZZ[e]/(e^n)$ is  a free abelian group of rank $n$ and, therefore, the standard representation of its Lie algebra of endomorphisms. 
The general case, for $r<\infty$  and arbitrary $n\geq r$,  follows from noticing that all $A\in \gln$ induce an even derivation $\delta(A)$ on $\wMn$:
\be
\left\{\begin{matrix}\delta(A)\bfu&:=&\hskip-105pt A\bfu,&&\forall \bfu\in M_n\cr\cr
\delta(A)(\bfv\w \bfw)&=&\delta(A)\bfv\w \bfw+\bfv\w \delta(A)\bfw,&&\forall \bfv,\bfw\in\wMn.
\end{matrix}\right.\label{eq0:2}
\ee
Since $\delta([A,B])=[\delta(A),\delta(B)]$, the map $A\mapsto \delta(A)_{|\bw^rM}$ makes $\bw^rM_n$ into an (irreducible) representation of  $\gln$. 
Then $B_{r,n}$ gets equipped with a $\gln$-module structure as well, due to the $\ZZ$-module isomorphism $B_{r,n}\sra \bw^rM_n$. Recall that the latter is the composition of the  Poincar\'e isomorphism, mapping  $B_{r,n}$ onto its singular homology $H_*(G(r,n),\ZZ)$, with the natural isomorphism  $H_*(G(r,n),\ZZ)\sra \bw^rH_*(\PP^{n-1},\ZZ)\cong \bw^rM_n$, as in \cite{SCHSD}, or \cite[diagramme (5.27)]{GatSal4}. 

For $r=\infty$, the fact that $B$ is a $\glin$-module is well known, and is due to the  isomorphism of $B$ with each degree of the fermionic Fock space which, roughly speaking, plays the role of an infinite exterior power.
 Its structure has been explicitly described by Date, Jimbo, Kashiwara and Miwa (DJKM) in \cite{DJKM01}, see also \cite[Formula (1.17)]{jimbomiwa} and \cite[p.~53]{KacRaRoz}, by computing the shape of a generating function $\Bcal(z,w)$ encoding the multiplication of any polynomial by elementary matrices $\Bcal_{ij}$ of infinite sizes.
 
In our contribution \cite{GatSal6} we determine the shape of the same generating function in the case $r<\infty$,  by using the formalism of Schubert derivations  in the sense of \cite{GatSal4}. The formula we obtain is  new (as far as we know), and has a classical flavor (occurring as a $2$-parameter deformation of the  Schur determinant occurring in Giambelli's formula). We then felt the need to show that our 
methods also work in the known case $r=\infty$. The output is the present 
paper, in which we offer an alternative deduction of the DJKM bosonic vertex representation of $\glin$, based on the extension of the Schubert derivations to an infinite wedge power. Our method to compute the $\gln$-structure of $B_{r,n}$ then works uniformly for all pairs $r\leq n$ ranging over $\NN\cup\{\infty\}$.
\claim{\bf Outline.} The main tool used in this paper is the notion of {\em Hasse--Schmidt} (HS) {\em derivation on an exterior algebra},  quickly recalled in Section \ref{sec:sec2}.
  Let $M:=\bigoplus_{i\in\ZZ}\ZZ\cdot b_i$ be a free abelian group with basis $
  \bfb:=(b_i)_{i\in\ZZ}$. A map $\Dcal(z):\wM\sra \wM[[z]]$ is said to be a 
  HS derivation on $\wM$ if $\Dcal(z)(u\w v)=\Dcal(z)u\w \Dcal(z)v$. In this paper, 
  we shall be  concerned mainly with the {\em Schubert derivations} (Section~
  \ref{sec:sec3}). They are denoted by $\sigma_+(z), \ovsig_+(z)$ and by $\sigma_-(z), \ovsig_-(z)$, where $\sigma_{\pm}(z):=\sum_{i\geq 0}\ovsig_{\pm i}z^{\pm i}\in \End_\ZZ(\wM)[[z^{\pm 1}]]$ are the unique HS derivations   such 
  that $\sigma_jb_i=b_{i+j}$, for all $i,j\in\ZZ$, and  $\ovsig_{\pm 1}(z)$ are their inverse in $\End(\bw M)[[z^{\pm 1}]]$. 
  
  The reason it is appropriate to call $\sigma_\pm(z)$ and $\ovsig_{\pm}(z)$
  Schubert derivations is explained in 
  \cite{GatSal4,pluckercone}. As a matter of fact, the operator $\sigma_i$ acting on $\bw^rM$, obeys the same 
  combinatorics enjoyed by the special Schubert cocycles in the cohomology ring of 
  a Grassmannian $G(r,n)$, for $r$ and $n$ big enough.
  
The Fermionic Fock space (FFS) come into the game in Section~\ref{sec:sec4}, playing the role of 
  something like $\bw^\infty M$ (often denoted in the literature by $
  \bw^{\infty/2}M$, to signify that is generated by semi-infinite exterior monomials, see e.g. \cite[Section 3]{KazLan01} or \cite[Section 1]{BlOko}). It is a notion for which there are excellent classical references in 
  the literature, such as \cite{Frenkel-zvi, Kacbeg, KacRaRoz}. However, to keep
  the exposition as self contained as possible, Section \ref{sec:sec4} supplies
  an alternative ad 
  hoc  algebraic construction of it, which widely suffices for our purposes and, possibly, may be useful for pedagogical ones.

 The extension of the Schubert derivations to the FFS is not entirely trivial, although not difficult. It turns out that the vertex operators occurring in the classical presentation of the 
Boson--Fermion correspondence can all be recovered by multiplying the four Schubert 
derivations $\sigma_\pm(z)$ and $\ovsig_\pm(z)$. For instance $\sigma_+(z)\ovsig_-(z)$ and $\ovsig_+(z)\sigma_-(z)$ are 
  basically the bosonic vertex operators acting on the Fock representation of the Heisenberg Lie algebra, \cite{Frenkel-zvi, KacRaRoz}, and  are described in 
Sections \ref{sec:sec6}--\ref{sec:sec7}. 
Let  $\delta(z,w):=\sum_{i,j\in\ZZ}\delta(\Bcal_{ij})z^iw^{-j}$, where $\Bcal_{ij}$ are as in the first part of this introduction. The main result of this paper is that the action of
  $
  \delta(z,w)
  $ on each degree of the FFS, 
  is proportional to the product 
  $$
  \sigma_+(z)\ovsig_-(z)\ovsig_+(w)\sigma_-(w)
  $$ of the four Schubert derivations.
The product $\sigma_-(z)\sigma_+
(w)$ commutes up to a rational factor, determined in Section \ref{sec:sec8} in order to achieve the final DJKM expression of Section~\ref{sec:sec9}, which is defined over  the integers.  Tensoring by $\QQ$, reading the expression as acting on $B\otimes_\ZZ\QQ$,  via its isomorphism with each degree of the FFS, and by  essentially  the same arguments as \cite[Theorem 7.7]{pluckercone}, one can check that: 
 $$
 \sigma_+(z)\ovsig_-(z)\ovsig_+(w)\sigma_-(w)=z \left(1-\displaystyle{w\over z}\right)^{-1}\sigma_+(z)\ovsig_+(w)\ovsig_-(z)\sigma_-(w)=z \left(1-\displaystyle{w\over z}\right)^{-1}\,\,\Gamma(z,w),
 $$ 
 where
\be
 \Gamma(z,w)=\exp\left(\sum_{i\geq 1}x_i(z^i-w^{i})\right)\exp\left(-\sum{1\over i}\left({1\over z^i}-{1\over w^i}\right){\partial\over \partial x_i}\right),\label{eq0:clasvo}
\ee
and the sequence $(x_1,x_2,\ldots)$ is defined through the equality $(1-e_1z+
e_2z^2-\cdots)\exp(\sum_i{x_iz^i})=1$. Then formula~(\ref{eq0:clasvo}) is precisely \cite[equation (5.33)]{KacRaRoz}. 

\claim{} We should finally remark that many of the tools employed in this paper within the 
framework of Schubert derivations have already been reviewed  
in other contributions (e.g. \cite{SCHSD,GatSal4, pluckercone, SCGA, GatSche03}), 
which we might well refer to. However, since the vocabulary of HS--derivations is 
not yet standard,   it seems motivated to recall the basic notions and
  facts  without saying anything more about proofs or the self-containedness of this first draft.
    
\medskip
\noindent{\bf Acknowledgments.} This work started during the stay of the first 
author at the Department of Mathematics of UNESP, S\~ao Jos\'e do Rio Preto, under 
the auspices of FAPESP, Processo n. 2016/03161--3. It continued under the partial 
support of INDAM-GNSAGA and the PRIN ``Geometria delle Variet\`a Algebriche''. A 
short visit of the first author to the second one was supported by the program 
``Finanziamento Diffuso della Ricerca'' of Politecnico di Torino. All these 
institutions are warmly acknowledged. For discussions and criticisms we want to 
thank primarily Israel Vainsencher, to whom this paper is dedicated on the occasion 
of his seventieth birthday,  as well as, in alphabetical order,  Carolina Araujo, Simon G.~Chiossi, Andr\'e Contiero,
Abramo Hefez, Marcos Jardim, Joachim Kock, Daniel Levcovitz, Simone Marchesi, Igor Mencattini, Piotr Pragacz, Andrea T. Ricolfi, Inna Scherbak and Aron Simis.

\section{Notation}\label{sec:nota}

\claim{}   A partition is a monotone non--increasing sequence $\blamb$ of non-negative integers $\lambda_1\geq\lambda_2\geq \ldots$ such that all terms are zero but finitely many. We denote  by $\ell(\blamb):=\sharp\{i\,|\,\lambda_i\neq 0\}$ its {\em length}. We denote by $\Pcal$ the set of all partitions and by $\Pcal_r$ the set of all partitions of length at most $r$. The partitions form an additive semigroup: if $\blamb,\bmu\in\Pcal$, then $\blamb+\bmu\in\Pcal$.
If $\blamb:=(\lambda_1,\lambda_2,\ldots)$, we denote by $\blamb^{(i)}$ the partition obtained by removing the $i$-th part:
$$
\blamb^{(i)}:=(\lambda_1\geq\lambda_{i-1}\geq\widehat{\lambda}_i\geq\lambda_{i+1}\geq\ldots,
$$
where \,\, $\widehat{}$\,\, means removed. 
By $(1^j)$ we mean the partition with $j$ parts equal to $1$.

\claim{} We denote by  $\wM=\bigoplus_{r\geq 0} \bw^rM$ the exterior algebra of a free abelian group $M:=\bigoplus_{i\in\ZZ}\ZZ\cdot b_i$ with basis $\bfb:=(b_i)_{i\in\ZZ}$. A typical element of $\bw^rM$ is a finite linear combination of monomials of the form $$
b_{i_{-1}}\w\cdots\w b_{i_{-r}}
$$
 with $\infty>i_{-1}>\cdots>i_{-r}>-\infty$.
 Given  $(m,r,\blamb)\in \ZZ\times \NN\times \Pcal_r$,  the following notation will be used:
\be 
\bfb^r_{m+\blamb}=b_{m+\lambda_1}\w b_{m-1+\lambda_2}\w\cdots\w b_{m-r+1+\lambda_r}\in \bw^rM_{\geq m-r+1}\subseteq \wM.\label{eq:notwbr}
\ee
\claim{} We  denote by $M_{\geq j}$ the sub-module of $M$ spanned by all $b_k$ with $k\geq j$, i.e. $M_{\geq j}:=\bigoplus_{i\geq j}\ZZ\cdot b_i$. In this case
$$
\bw^rM_{\geq j}:=\bigoplus_{\blamb\in\Pcal_r}\ZZ\bfb^r_{r-1+j+\blamb}\qquad \mathrm{and}\qquad  \bw M_{\geq j}=\bigoplus_{r\geq 0}\bw^rM_{\geq j}.
$$

\section{Hasse-Schmidt Derivations on Exterior Algebras}\label{sec:sec2}

Main detailed references for this  section are  \cite{SCHSD,GatSal4,pluckercone}.

\claim{}   A {\em Hasse-Schmidt} (HS) {\em derivation} on $\wM$ is a $\ZZ$-linear map  $\Dcal(z):\wM\sra \wM[[z]]$, such that for all $\bfu,\bfv\in\wM$:
\be
\Dcal(z)(\bfu\w \bfv)=\Dcal(z)\bfu\w\Dcal(z)\bfv.\label{eq:HSder}
\ee
If $D_i\in\End_\ZZ(\wM)$  is such that $\Dcal(z)=\sum_{i\geq 0}D_iz^i$,  then  equation \eqref{eq:HSder} is equivalent 
to 
\be 
D_i(\bfu\w \bfv)=\sum_{j=0}^i\Dcal_i\bfu\w D_{i-j}\bfv.
\ee
If $D_0$ is invertible, up to termwise multiplying $\Dcal(z)$ by $D_0^{-1}$, we may assume that $D_0=\mathrm{id}_{\wM}$. Thus $\Dcal(z)$ is invertible in $\End_\ZZ(\wM)[[z]]$. An easy check shows that the formal inverse $\ovDcal(z):=\sum_{j\geq 0}(-1)^j\ovDcal_jz^j$ is a HS-derivation as well.  

%The relation
%\be \ovDcal(z)\Dcal(z)=\Dcal(z)\ovDcal(z)=1,\label{eq:HSinv}
%\ee
%is equivalent to the following set of equalities
%\be
%\sum_{i=0}^j(-1)^i\Dcal_j\ovDcal_{j-i}=0.
%\ee
A main tool of this paper is:
\bclm{\bf Proposition.} {\em The integration by parts formulas hold:
\begin{eqnarray}
\Dcal(z)\bfu\w \bfv&=&\Dcal(z)(\bfu\w \ovDcal(z)\bfv),\label{eq:intprt1}\\ \cr
\bfu\w \ovDcal(z)\bfv&=&\ovDcal(z)(\Dcal(z)\bfu\w \bfv).
\end{eqnarray}
}
\eclm
\proof  Straightforward from the definition.\qed

\claim{\bf Duality.}  Let $\beta_j: M\sra \ZZ$ be the unique linear form such that $\beta_j(b_i)=\delta_{ij}$. The {\em restricted dual} of $M$ is  $M^*:=\bigoplus_{j\in\ZZ}\ZZ\cdot\beta_j$. Recall the natural identification between $\bw^rM^*$ and $(\bw^rM)^*$:
$$
\beta_{i_1}\w\cdots\w\beta_{i_r}(b_{j_1}\w\cdots\w b_{j_r})=\left|\begin{matrix}\beta_{i_1}(b_{j_1})&\cdots&\beta_{i_1}(b_{j_r})\cr
\vdots&\ddots&\vdots\cr
\beta_{i_r}(b_{j_1})&\cdots&\beta_{i_r}(b_{j_r})
\end{matrix}\right|.
$$
The {\em contraction}   of $u\in \bw^rM$ against $\beta\in M^*$ is the unique vector $\beta\lrcorner \bfu\in \bw^{r-1}M$ such that  the equality
$$
\eta(\beta\lrcorner u)=(\beta\w \eta)(u),
$$
holds for all $\eta\in\bw^{r-1}M^*$.
\bclm{\bf Definition.} {\em 
The map $
\Dcal^T(z)=\sum_{i\geq 0}\Dcal_i^Tz^i:\wM^*\sra \wM^*[[z]]
$ 
such that
$$
(\Dcal^T(z)\eta)(\bfu)=\eta (\Dcal(z)\bfu),
$$
is called the {\em transpose} of the HS-derivation $\Dcal(z)$.
 }
 \eclm
By~\cite[Proposition 2.8]{pluckercone}, it follows that $\Dcal^T(z)(\eta_1\w\eta_2)=\Dcal^T(z)\eta_1\w\Dcal^T(z)(\eta_2)$, for all $\eta_1,\eta_2\in \bw M^*$, i.e that $\Dcal^T(z)$ is a HS derivation on $\bw M^*$.
%\claim{} An easy argument, see e.g., shows that there is one and only one HS-derivation $\Dcal^f(z):\wM\sra\End_\ZZ(\wM)[[z]]$  extending an arbitrary $f(z)\in \End_\ZZ(\wM)[[z]]$. On homogeneus elements is given by 
%\be 
%\Dcal^f(z)(u_1\w\cdots\w u_r)=f(z)u_1\w\cdots\w f(z)u_r
%\ee
%for all $(u_1,\ldots, u_r)\in M^r$.

\section{Schubert Derivations}\label{sec:sec3}

\bclm{\bf Definition.}\label{def:schubder}{\em  The {\em Schubert derivations} are the  unique HS-derivations
\be
\sigma_+(z), \sigma_-(z):\wM\sra \wM[[z^{\pm 1}]]\label{eq:schubder}
\ee
such that $\sigma_\pm(z)b_j=\sum_{i\geq 0}b_{j\pm i}z^{\pm i}$ and their formal inverses $\ovsig_{\pm}(z)$ in $\End(\bw M)[[z^{\pm 1}]]$.
Their formal inverses $\ovsig_\pm(z)\in \End(\wM)[[z]]$ are the unique HS-derivations $\wM \sra \wM[[z]]$  such that 
\be
\ovsig_\pm(z)b_j=b_j-b_{j\pm 1}z.\label{eq:invsder}
\ee
}
\eclm
\claim{\bf Remark.} To save notation, we preferred to write  $\sigma_-(z),\ovsig_-(z)$  rather than the more precise $\sigma_-(z^{-1})$ and $\ovsig_-(z^{-1})$,  hoping that the subscript ``${_-}$'' to $\sigma$ may  help to avoid possible confusions.

 Put $\sigma_\pm(z)=\sum_{j\geq 0}\sigma_{\pm j}z^{\pm j}$ and $\ovsig_\pm(z)=\sum_{j\geq 0}(-1)^j\ovsig_{\pm j}z^{\pm j}$. Then:
\be
\sigma_{ i}b_j=b_{i+j},\qquad\forall i,j\in\ZZ,\label{eq3:32}
\ee
while $\ovsig_{i}u=0$ if $u\in \bw^{\leq\, \left|i\right|-1}M$, for all $i\in\ZZ$
(Cf. \cite[Secs. 3.1--3.2]{pluckercone}).
\bclm{\bf Remark.} The operator $\sigma_i$ defined on $M$ are precisely the shift 
operators $\Lambda_i$ as in \cite[p.~32]{KacRaRoz}. The only difference is that i) 
we extend them to all the exterior algebra of $M$ (and then to the associated 
fermionic Fock space) embedding them into a Schubert derivation; ii) due to i), we 
preferred to use the notation $\sigma_i$ to emphasize the interpretation in terms 
of Schubert calculus. The shift operators $\sigma_i$ acts on $b_0$ as the cap product of  the class of a linear space of codimension $i$ with the fundamental class of some $\PP^n$ (being $0$ if $i>n$).

\eclm

\bclm{\bf Remark.}\label{rmk:rmk1} Let $\sigma^*_+(w):\wM^*\sra\wM^*[[w]]$ be the Schubert derivation on $\wM^*$, i.e. the unique HS-derivation such that
$$
\sigma^*(w)\beta_j=\sum_{i\geq 0}\sigma_i^*\beta_j\cdot w^i=\sum_{i\geq 0}\beta_{j+i}w^i.
$$ Its inverse $\ovsig^*_+(w)$  is the unique HS-derivation on $\wM^*$ such that 
$$
\ovsig_+^*(w)\beta_j=\sum_{i\geq 0}(-1)^i\ovsig_i^*\beta_jw^i=\beta_j-\beta_{j+1}w.
$$
An easy check shows that $\sigma_-(z)=\sum \sigma_{-i}z^{-i}=\sigma_+^{*T}
(w)_{|{w=z^{-1}}}$. Similarly $\ovsig_-(z)={\ovsig_+^{*T}(w)}_{|{w=z^{-1}}}:\wM \sra ù
\wM[[z^{-1}]]$.
\eclm

%\bclm{\bf Definition.} {\em 
%Let $\sigma_-(\z)=\sum \sigma_{-i}z^{-i}:=\sigma_+^{*T}(z):\wM \sra \wM[[z^{-1}]]$ and $\ovsig_-(\z)=\sum \ovsig_{-i}z^{-i}:=\ovsig_+^{*T}(z^{-1}):\wM \sra \wM[[z^{-1}]]$ 
%}
%\eclm
%
%Thus, $\sigma_-(z):=\sum_{i\geq 0} \sigma_{-i}z^{-i}$  is the unique HS-derivation of $\wM$ such that 
%$$
%\sigma_-(z)b_j=\sum_{i\geq 0}\sigma_{-i}b_{j}z^{-i}=\sum_{i\geq 0}b_{j-i}z^{-i}
%$$
%while $\ovsig_-(z):\wM\sra \wM[z^{-1}]$ is the unique HS derivation on $\wM$ such that
%$$
%\ovsig_-(z)b_j=\sum_{i\geq 0}\ovsig_{-i}b_{j}z^{-i}=b_j-b_{j-1}z^{-1}.
%$$
\bclm{\bf Proposition.} {\em The following equalities hold:
\begin{eqnarray}
\sigma^T_-(z)\beta_j&=&\sum_{i\geq 0}\beta_{j+i}z^{-i}={\sigma_+^*(w)\beta_j}_{|w=z^{-1}},
\\ \cr
\ovsig_-^T(z)\beta_j&=&\beta_j-\beta_{j+1}z^{-1}={\ovsig^*(w)\beta_j}_{|w=z^{-1}}.
\end{eqnarray}
}
\eclm
\proof It follows from the definition.\qed

\section{Fermionic Fock Space}\label{sec:sec4}

\claim{} There are several excellent references concerning the definition of  the fermionic Fock space \cite[Ch. 5]{Frenkel-zvi} or \cite{Kacbeg,KacRaRoz}. It amounts to the rigorous formalization of the idea of an infinite exterior power.   We propose here an elementary algebraic construction of it,  that suffices for our purposes.

Let $[M]$ be the free $\ZZ$-module generated by the basis $\wb:=(\wb_m)_{m\in\ZZ}$. Identify $[M]$ with a sub-module of the tensor product $\wM\otimes_\ZZ[M]$ via the map $\wb_m\mapsto 1\otimes \wb_m$. Let $W$ be the $\wM$--submodule  of $\wM\otimes_\ZZ[M]$ generated by all the expressions $\{b_m\otimes \wb_{m-1}-\wb_m, b_m\otimes [b]_m\}_{m\in\ZZ}$. In formulas:
\[
W:=\wM\otimes \big(b_m\otimes \wb_{m-1}-\wb_m\big)+ \wM\otimes\big( b_m\otimes [b]_m\big).
\]
\bclm{\bf Definition.} {\em The fermionic Fock space is the $\wM$-module
\be
F:=F(M):={\wM\otimes_\ZZ[M]\over W}.
\ee
}
\eclm
Let $\wM\otimes_\ZZ[M]\sra F$ be the canonical projection. The class  of $u\otimes \wb_m$ in $F$ will be denoted  $u\w [b]_m$. Thus the equalities $b_m\w \wb_m=0$ and $b_m\w\wb_{m-1}=\wb_m$ hold in $F$. For all $m\in \ZZ$ and $\blamb\in\Pcal$ let, by definition
$$
\wb_{m+\blamb}:=\bfb^r_{m+\blamb}\w \wb_{m-r},
$$
where $r$ is any positive integer such that $\ell(\blamb)\leq r$. Then $F$ is a graded $\wM$-module:
$$
F:=\bigoplus_{m\in\ZZ}F_m,
$$
where
\be 
F_m:=\bigoplus_{\blamb\in\Pcal}\ZZ\wb_{m+\blamb}=\bigoplus_{r\geq 0}\bigoplus_{\blamb\in\Pcal_r}\ZZ\bfb^r_{m+\blamb}\w \wb_{m-r}.\label{eq:Fm}
\ee
\bclm{\bf Definition.} {\em The {\em fermionic Fock space} of charge $m$ is the module $F_m$ as in~\eqref{eq:Fm}}\cite[p. 36]{KacRaRoz}.
\eclm
\bclm{\bf Proposition.} {\em The equality  $b_j\w\wb_m=0$ holds for all $j\leq m$.
}
\eclm
\proof Indeed:

\smallskip
$
b_j\w\wb_m =b_j\w b_m\w\cdots\w b_j\w\wb_{j-1}=\pm b_j\w b_j\w b_m\w\cdots\w \wb_{j-1}=0.
$\qed

\bclm{\bf Proposition.} {\em The image of the map $\bw^rM\otimes F_m\sra F$ given by $(u,v)\mapsto u\w v$ is contained in $F_{m+r}$.}
\eclm
\proof
Let $b_{i_1}\w\cdots\w  b_{i_r}\in \bw^rM$ with $i_1>\cdots>i_r$. Write 
$$
i_1=m+1+r+\lambda_1,\quad \cdots,\quad i_r=m+1+\lambda_r.
$$
If $\lambda_1\geq\cdots\geq \lambda_r$, then
$$
b_{i_1}\w\cdots\w  b_{i_r}\w[b]_m=\bfb^r_{m+1+r+\blamb}\w \wb_m\in F_{m+r},
$$
otherwise the product is zero.
\qed

\claim{} In particular each $\bfu\in M$ defines an action $\bfu\,\w:F_m\sra F_{m+1}$ given by
\be
\bfb^r_{m+\blamb}\w [b]_{m-r}\longmapsto (u\w \bfb^r_{m+\blamb})\w [b]_{m-r}.\label{eq4:wedging}
\ee

%Let $\Pcal$ be the set of all partitions and  $\Pcal_r$ the set of all partitions of length at most $r$. Denote by $M_{\geq k}:=\bigoplus_{i\geq k} \ZZ\cdot b_i$. If $\blamb\in\Pcal_r$ let
%\be 
%\bfb^r_{m+\blamb}:=b_{m+\lambda_1}\w b_{m-1+\lambda_2}\w\cdots\w b_{m-r+1+\lambda_r}\in \bw^rM_{\geq m-r+1}
%\ee
%Let us now denote by $\Bcal$ the free $\ZZ$-module generated by all the expressions $\{\wb_m\,|\, m\in\ZZ\}$.  We let $F$ to be the quotient of the tensor product $\wM\otimes \Bcal$ modulo the relations $b_m\otimes \wb_{m-1}=\wb_m$ and $b_j\w\wb_m=0$ for all $j\leq m$. The class of $u\otimes \wb_m$, for all $u\in\wM$ will be denoted by $u\w \wb_m$.
%If $\blamb\in\Pcal_r$, let
%$$
%\wb_{m+\blamb}=\bfb^r_{m+\blamb}\w\wb_{m-r}
%$$
%Notice that due to the relations the expression does not depend on the choice of the integer $r$ provided it is greater than $\ell(\blamb)$.
%Notice that due to the relations $\wb_m$ can be expressed as an infinite wedge monomial:
%$$
%\wb_m=b_m\w \wb_{m-1}=b_m\w b_{m-1}\w\wb_{m-2}=b_m\w b_{m-1}\w b_{m-2}\w\cdots
%$$
%Let
%$$
%F_m:=\bigoplus\ZZ \wb_{m+\blamb}
%$$
%Then,  clearly
%$$
%F:=\bigoplus_{m\in\ZZ}F_m
%$$
%and $F_m$ is said to be the fermionic Fock space of charge $m$. The element $\wb_m$ is usually denoted in literature by the symbol $|m\rangle$, the {\em vacuum} in charge $m$.
Similarly,  one may consider a contraction  action of $\bw M^*$ on $F$, mapping $F_m\mapsto F_{m-1}$. Define the contraction of $\wb_{m+\blamb}$ against $\beta_j$ as follows: choose $r$ such that $\ell(\blamb)\leq r$ and $m-r\leq j$. Declare that
\be
\beta_j\lrcorner \wb_{m+\blamb}=\beta_j\lrcorner \bfb^r_{m+\blamb}\w \wb_{m-r},\label{eq4:contra}
\ee
and extend by linearity. The definition does not depend on the choice of $r>\max(\ell(\blamb), m-j)$.
For instance
$$
\beta_m\lrcorner \wb_m=\beta_m\lrcorner (b_m\w \wb_{m-1})=\wb_{m-1}.
$$

%Let $\Ccal:=\Ccal(M\oplus M^*)$ be the Clifford algebra associated to the canonical inner product 
%$
%<u\oplus \alpha,v\oplus\beta>=\alpha(v)+\beta(u).
%$
%By identifying $b_i$ with $b_i\oplus {\mathbf 0}$ and $\beta_j$ with ${\mathbf 0}\oplus \beta_j$, it follows that $(\bfb,{\bm \beta}):=(b_i,\beta_j)_{i,j\in\ZZ}$ is a $<,>$-orthonormal basis of $M\oplus M^*$.
%\bclm{\bf Proposition.} {\em The maps $\beta_j\wb_{m+\blamb}=\beta_j\lrcorner \wb_{m+\blamb}$ and $b_i\wb_{m+\blamb}=b_i\w \wb_{m+\blamb}$ define an action of the Clifford algebra $\Ccal$ on $F$.
%}
%\eclm
%\proof
%It suffices to show that for all $i,j,m\in\ZZ$ and $\blamb\in\Pcal_r$
%$$
%\beta_j\lrcorner (b_i\w \wb_{m+\blamb})+b_i\w (\beta_j\lrcorner\wb_{m+\blamb})=\delta_{i,j}\wb_{m+\blamb}
%$$
%Choose $r$ big enough such that $\ell(\blamb)\leq r$ and $m-r-1\leq j$, so that, by definition:
%\begin{eqnarray*}
%\beta_j\lrcorner (b_i\w \wb_{m+\blamb})+b_i\w (\beta_j\lrcorner\wb_{m+\blamb})&=&[\beta_j\lrcorner (b_i\w \bfb^r_{m+\blamb})\cr &+&b_i\w (\beta_j\lrcorner \wb_{m+\blamb}^r)]\w \wb_{m-r}\cr\cr
%&=&\delta_{i,j}\wb^r_{m+\blamb}\w \wb_{m-r}=\delta_{i,j}\wb_{m+\blamb}.\hskip 35pt \qed
%\end{eqnarray*}

\claim{\bf Remark.} \label{rmkfermion} Using the {\em wedging} and {\em contraction} operators \eqref{eq4:wedging} and \eqref{eq4:contra},  is easy to show that $F$ is an irreducible representation of a canonical Clifford algebra on $M\oplus M^*$, called Fock fermionic representation in \cite[Section 5.]{Frenkel-zvi}, which motivates the terminology we adopted.

\claim{\bf Duality.} \label{sec:duality} We denote by $F^*$ the restricted dual of $F$, contructed out of the restricted dual of $M^*$, precisely as one did for $F$.  The typical element of $F$ is of the form
$$
\wbet_{m+\blamb}=\wbet^r_{m+\blamb}\w \wbet_{m-r}.
$$
The duality pairing $F^*\times F\sra \ZZ$ is defined by
$
\wbet_{m+\blamb}(\wb_{n+\bmu})=\delta_{m,n}\delta_{\blamb,\bmu}
$.
It extends the natural duality between $\bw^rM$ and $\bw^rM^*$. The contraction of $f\in F_m$ against $\beta\in M^*$ is defined by
$$
\eta(\beta\lrcorner f)=(\beta\w \eta)(f),
$$
for all $\eta\in F^*_{m-1}$.
More explicitly, keeping into account that $f$ is a finite sum of elements of the form $\wb_{m+\blamb}$, one has:
$$
\beta\lrcorner \wb_{m+\blamb}=\beta\lrcorner\bfb^r_{m+\blamb}\w \wb_{m-r}+(-1)^{r-1}\bfb^r_{m+\blamb}\w (\beta\lrcorner \wb_{-m-r}).
$$
Since $\beta$ is a finite sum $\sum a_i\beta_i$, the contraction $\beta\lrcorner \wb_{-m-r}$ is a finite sum.
%\claim{\bf Boson-Fermion correspondence.} Let $B:=\ZZ[e_1,e_2,\ldots]$ be the polynomial ring in infinitely many indeterminates $(e_1,e_2,\ldots)$. Denote $E(z)=1-e_1z+e_2z^2+\cdots$ and define $(h_j)_{j\in\ZZ}$ via the relation
%\be 
%\sum_{j\in\ZZ}h_jz^j:={1\over E(z)}
%\ee
%Then $B:=\ZZ[h_1,h_2,\ldots]$ as well. Moreover 
%$$
%B=\bigoplus \ZZ\Delta_\blamb(H)
%$$
%where $\Delta_\blamb(H)=\det(h_{\lambda_j+i-j})_{1\leq i,j\leq \ell(\blamb)}$. The unique $\ZZ$-module isomorphism
%\be 
%B[u,u^{-1}]\lra \bw^{\infty/2}M
%\ee
%mapping $u^m\Delta_\blamb(H)\mapsto \wb_{m+\lambda}$ is called the {\em Boson-Fermion correspondence}. In fact $F$ is a $B[u,u^{-1}]$ module generated by $\wb_0$, by setting
%$$
%u^m\Delta_\blamb(H)\wb_{n+\bmu}=\Delta_\blamb(H)\Delta_\bmu(H)\wb_{m+n}=u^{m+n}\Delta_\blamb(H)\Delta_\bmu(H)\wb_0
%$$
\section{Extending HS-derivations to $F$}\label{sec:sec5}
 This section is devoted   to extend the Schubert derivations \eqref{def:schubder}, and their transposed, to suitable maps  $F\sra F[[z]]$ and $F^*\sra F^*[[z^{-1}]]$. The purpose is to (re)-discover the {\em bosonic} vertex operators as in \cite[Theorem 5.1]{KacRaRoz} or \cite[p.~92]{Kacbeg}. This will supply  an alternative way to look at the bosonic vertex representation of the Lie algebra $gl_\infty(\ZZ)$, due to Date-Jimbo-Kashiwara-Miwa \cite{DJKM01,jimbomiwa}. Although those authors worked over the complex numbers, we work over the integers because it is sufficient for our purposes. 
 
 The sought for extension  of the Schubert derivations to $F$ will be attained looking at each degree $F_m$ at a time.
\claim{}\label{subs:s+ovs+}
Let $\sigma_+(z), \sigma_-(z):\wM\sra \wM[[z^{\pm 1}]]$ be the Schubert derivations and  $\ovsig_+(z), \ovsig_-(z):\wM\sra \wM[[z^{\pm 1}]]$ their inverses as in Section~\ref{def:schubder}. We first extend them to $\ZZ$-basis elements of $[M]$ and then we extend to all $F$ by  mimicking the typical behavior of an  algebra homomorphism.

\bclm{} \label{def:preextend} {Let us begin to extend the definition of $\sigma_{\pm}(z)$ and $\ovsig_{\pm}(z)$ to elements of $[M]$ as follows:
\be
\sigma_+(z)\wb_m:=(\sigma_+(z)b_m)\w \wb_{m-1},
\ee
and 
\be
\ovsig_+(z)\wb_m:=\sum_{j\geq 0}(-1)^j\wb_{m+(1^j)}z^j\label{eq:defon1M}
\ee
for all $m\in\ZZ$.
We demand, instead, that  $\sigma_-(z)$ and $\ovsig_-(z)$ act on $[M]$  as the identity:
\be 
\sigma_-(z)\wb_m=\wb_m\qquad \mathrm{and}\qquad \ovsig_-(z)\wb_m=\wb_m,\label{eq5:sigact}
\ee
for all $m\in\ZZ$.
}
\eclm
Notice that \eqref{eq:defon1M} can be equivalently written in $F_m$ as:
$$
\sigma_+(z)\wb_m=\sum_{j\geq 0}(-1)^j\wb^j_{m+(1^j)}\w \wb_{m-j}z^j.
$$
We now extend the Schubert derivations, as in \ref{def:schubder},  to all $F_m$.

\bclm{\bf Definition.}\label{def:extend}  {\em The extension of the Schubert derivations  $\sigma_\pm(z),\ovsig_\pm(z)$  to $\ZZ$-linear maps  $F_m\sra F_m[[z^{\pm 1}]]$ is defined by:
\begin{eqnarray}
\sigma_\pm(z)\wb_{m+\blamb}&=&\sigma_\pm(z)(\bfb_{m+\blamb}^r\w\wb_{m-r})\cr\cr &=&\sigma_\pm(z)(\bfb_{m+\blamb}^r)\w\sigma_\pm(z)(\wb_{m-r}),\label{eq:defs+}
\end{eqnarray}
and
\begin{eqnarray}
\ovsig_\pm(z)(\wb_{m+\blamb})&=&\ovsig_\pm(z)(\wb_{m+\blamb}^r\w\ovsig_{\pm}(z)\wb_{m-r})\cr\cr &=&\ovsig_+(z)\wb^{r}_{m+\blamb}\w \ovsig_+(z)\wb_{m-r}.\label{eq:defovs+}\end{eqnarray}
}
\eclm
%\bclm{\bf Proposition.} {\em For all $m\in\ZZ$,  $\sigma_+(z)\wb_m=\sigma_+(z)b_m\w \wb_{m-1}$.
%}
%\eclm
%\proof In fact, by \eqref{eq:defs+}
%\begin{eqnarray}
%\sigma_+(z)\wb_m&=&\sigma_+(z)(b_m\w b_{m-1})\w \wb_{m-2}\cr\cr
%&=&\sigma_+(z)b_m\w (b_{m-1}+z\ovsig_+(z)b_m)\w \wb_{m-2}\cr\cr
%&=&\sigma_+(z)b_m\w b_{m-1}\w \wb_{m-2}=\sigma_+(z)b_m\w\wb_{m-1}
%\end{eqnarray}
%\qed
\bclm{\bf Proposition.}\label{prop:corsibm}
{\em For all $i\in\ZZ$, $\sigma_i\wb_m=b_{m+i}\w \wb_{m-1}$ and is thence zero if $i<0$.}
\eclm
\proof
If $i\geq 0$, $\sigma_{i}\wb_m$ is the coefficient of $z^{i}$ in the expression
$$
\sigma_\pm(z)b_m\w \wb_{m-1}=\sum_{i\geq 0}\sigma_ib_mz^i\w\wb_{m-1}=\sum_{i\geq 0}b_{m+i}\w \wb_{m-1}z^i.
$$
If $-i>0$, instead, $\sigma_{-i}\wb_m$ is the coefficient of $z^{-i}$ in the right--hand side of the equation $\sigma_-(z)\wb_m=\wb_m$, which is zero as stated.\qed
\bclm{\bf Proposition.} {\em For all $i\geq 0$, the maps $\sigma_{\pm i}:F_m\sra F_m$ satisfy the  $i$-th order Leibniz rule
\be
\sigma_{\pm i}(\bfb^r_{m+\blamb}\w \wb_{m-r})=\sum_{j=0}^i\sigma_{\pm j}\bfb^r_{m+\blamb}\w \sigma_{\pm i\mp j}\wb_{m-r}.\label{eq:s_ileib}
\ee
}
\eclm
\proof
In fact, the left--hand side of \eqref{eq:s_ileib} is the coefficient of $z^{\pm i}$ of the expression
$
\sigma_\pm(z)(\bfb^r_{m+\blamb}\w \wb_{m-r})
$, 
which by definition is
$$
\sigma_\pm(z)(\bfb^r_{m+\blamb})\w \sigma_\pm(z)\wb_{m-r}.
$$
The equality
$$
\sigma_-(z)\bfb^r_{m+\blamb}\w \sigma_-(z)\wb_{m-r}=\sigma_-(z)\bfb^r_{m+\blamb}\w \wb_{m-r}
$$
implies 
$$
\sigma_{-i}\bfb^r_{m+\blamb}\w \wb_{m-r}=\sigma_{-i}\bfb^r_{m+\blamb}\w \wb_{m-r},
$$
as all the other terms involving $\ovsig_{-j}\wb_{m-r}$ vanish for $j>0$.
Moreover, using Proposition \ref{prop:corsibm}, one sees that the coefficient of $z^{ i}$ in $\sigma_+(z)(\bfb^r_{m+\blamb}\w \wb_{m-r})$ is
\begin{eqnarray}
\sigma_{i}(\bfb^r_{m+\blamb}\w \wb_{m-r}&=&\sum_{j=0}^i\sigma_j\bfb^r_{m+\blamb}\w \sigma_{i-j}\wb_{m-r}=\sum_{j=0}^i\sigma_j\bfb^r_{m+\blamb}\w \sigma_{i-j}\wb_{m-r}
\end{eqnarray}
as desired.\qed
\bclm{\bf Proposition.} {\em The following equalities hold:
$$
\ovsig_{\pm i}(\bfb^r_{m+\blamb}\w\wb_{m-r})=\sum_{j=0}^i\ovsig_{\pm j}\bfb^r_{m+\blamb}\w \ovsig_{\pm i\mp j}\wb_{m-r}
$$
}
\eclm
\proof The formula readily follows by comparing the coefficient of $z^{\pm i }$ on either side of formula \eqref{eq:defovs+}.\qed

\smallskip
\noindent
In particular $\ovsig_{-j}(\bfb^r_{m+\blamb}\w \wb_{m-r})=\ovsig_{-j}(\bfb^r_{m+\blamb})\w\wb_{m-r}$

\bclm{\bf Lemma.}\label{lem:le} {\em
Let  $\bfz_r:=(z_1,\ldots,z_r)$ be indeterminates over $\ZZ$. Then:
\be 
\ovsig_+(z_1)\cdots\ovsig_+(z_r)b_{m-r}=b_{m-r}+\sum_{j=1}^r(-1)^je_j(\bfz_r)b_{m-r+j},\label{eq:ovsz1zr}
\ee
where $e_i(\bfz_r)$ is the $i$-th elementary symmetric polynomial in $(z_1,\ldots,z_r)$.
}
\eclm
\proof If $r=1$, $\ovsig_+(z_1)b_{m-1}=b_{m-1}-b_mz_1=b_{m-1}-e_1(z_1)b_m$, showing that   Lemma~\ref{lem:le} holds in this case. Suppose that \eqref{eq:ovsz1zr}  holds for all $1\leq s\leq r-1$. Then
\begin{eqnarray*}
&&\ovsig_+(z_1)\ovsig_+(z_2)\cdots\ovsig_+(z_r)b_{m-r}\cr\cr &=&\ovsig_+(z_1)(b_{m-r}-e_1(z_2,\ldots, z_r)b_{m-r+1}+\cdots+(-1)^rz_2\cdots z_rb_{m-1})\cr\cr
&=&\sum_{j=0}(-1)^je_j(z_2,\ldots,z_r)(b_{m-r-j}-z_1b_{m-r-j+1})\cr
&=&b_{m-r}+\sum_{j=1}^r(-1)^je_j(\bfz_r)b_{m-r+j},
\end{eqnarray*}
as desired. \qed
\bclm{\bf Proposition.}\label{cor:cor66} {\em For all $m\in\ZZ$ and $r\geq 0$:
\be
\bfb^r_m\w \ovsig_+(z_1)\cdots\ovsig_+(z_r)\wb_{m-r}=\wb_m.\label{eq:bmws+bm}
\ee
}
\eclm
\proof
One has
$$
b_m\w\cdots\w b_{m-r+1}\w  \sum_{j=0}^r(-1)^j e_j(\bfz_r)b_{m-r+j} \w \sum_{j=0}^r (-1)^je_j(\bfz_r)b_{m-r-1+j}\w\cdots=\wb_m.
$$
\bclm{\bf Corollary.}\label{prop:bmsovsbm}
{\em For all $m\in\ZZ$:
$$
b_m\w\ovsig_+(z)\wb_{m-1}=\wb_m.
$$
}
\eclm
\proof
In fact the typical coefficient of $z^i$,  $i>0$, is the sum of monomials of the form $b_m\w b_m\w\cdots=0$, so that the only surviving summand is $b_m\w\wb_{m-1}=\wb_m$.\qed

\bclm{\bf Proposition.} {\em The maps $\ovsig_+(z),\sigma_+(z): F_m\sra F_m[[z]]$  and $\ovsig_-(z^{-1}),\sigma_-(z): F_m\sra F_m[[z^{-1}]]$ are one the inverse of the other.
} \eclm
\proof The formal power series $\sigma_+(z)$ and $\ovsig_+(z)$ are both invertible in $F_m[[z]]$, as the constant term is the identity of $F_m$. It then suffices to show that $\ovsig_+(z)$ is the left inverse of $\sigma_+(z)$, because in this case it must coincide with its inverse. Now for all $m\in\ZZ$ and each $\blamb\in\Pcal$, 
$$
\ovsig_+(z)\sigma_+(z)\wb_{m+\blamb}=\ovsig_+(z)\sigma_+(z)(\bfb^r_{m+\blamb}\w \wb_{m-r})=\bfb^r_{m+\blamb}\w \ovsig_+(z)\sigma_+(z)\wb_{m-r}
$$
for all $r\geq \ell(\blamb)$, where we used the fact that $\ovsig_+(z),\sigma_+(z)$ are HS--derivations on $\bw M_{\geq m-r+1}$, one inverse of the other.
It will then suffice to show that $\ovsig_+(z)\sigma_+(z)$ acts on $\wb_m$ as the identity, for all $m\geq 0$. Using Corollary \ref{prop:corsibm}, one easily gets:
$$
\ovsig_+(z)\sigma_+(z)\wb_m=\ovsig_+(z)\sigma_+(z)b_m\w \wb_{m-1}=b_m\w \ovsig_+(z)\wb_{m-1}=\wb_m 
$$
Similarly, since
$$
\ovsig_-(z)\sigma_-(z)b_j=b_j,
$$
for all $j\in \ZZ$, then
\begin{eqnarray*}
 \ovsig_-(z)\sigma_-(z)\wb_{m+\blamb}&=&\ovsig_-(z)[\sigma_-(z)\bfb]_{m+\blamb}\cr
 &=&[\ovsig_-(z)\sigma_-(z)\bfb]_{m+\blamb}=\wb_{m+\blamb}.\hskip205pt\qed
\end{eqnarray*}

\bclm{\bf Proposition.} {\em For all partitions $\blamb$ of length  at most $r$ and all $m\in\ZZ$, the integration by parts formulas hold:
\be
\sigma_+(z)\bfb^r_{m+\bmu}\w \wb_{m-r+\blamb}=\sigma_+(z)(\bfb^r_{m+\bmu}\w \ovsig_+(z)\wb_{m-r+\blamb});\label{eq:8}
\ee
\be
\bfb^r_{m+\bmu}\w \ovsig_+(z)\wb_{m-r+\blamb}=\ovsig_+(z)(\sigma_+(z)\bfb^r_{m+\bmu}\w \wb_{m-r+\blamb}).\label{eq:9}
\ee
}
\eclm
\proof Obvious.\qed
\bclm{\bf Proposition.}\label{prop:prop1} Let $i_1,\ldots,i_r$ be an index sequence of length $r\geq 1$. Then 
\be
\sigma_{i_1}\cdots \sigma_{i_r}\wb_m=(\sigma_{i_1}\cdots\sigma_{i_r}\bfb^r_m)\w \wb_{m-r}.\label{eq:sigmaI}
\ee
\eclm
\proof Let $(z_1,\ldots,z_r)$ be indeterminates over $\ZZ$. 
Then $\sigma_{i_1}\cdots\sigma_{i_r}\wb_m$ is the coefficient of $z_1^{i_1}\cdots z_r^{i_r}$ in the expansion of $\sigma_+(z_1)\cdots \sigma_+(z_1)\wb_m$. Now, by virtue of Corollary~\ref{cor:cor66}:
\begin{eqnarray}
\sigma_+(z_1)\cdots\sigma_+(z_r)\wb_m=\sigma_+(z_1)\cdots\sigma_+(z_r)(\bfb^r_m\w \ovsig_+(z_1)\cdots\ovsig_+(z_r)\wb_{m-r})
\end{eqnarray}
As $\sigma_+(z_i)$ is a HS-derivation, one obtains
\be
\sigma_+(z_1)\cdots\sigma_+(z_r)\wb_m=\sigma_+(z_1)\cdots\sigma_+(z_r)\bfb^r_m\w \wb_{m-r}.\label{eq:eitside}
\ee
Comparing the coefficient of $z_1^{i_1}\cdots z_r^{i_r}$ on either side of equality~\eqref{eq:eitside}, yields formula~\eqref{eq:sigmaI}.\qed

\bclm{\bf Proposition.} {\em For all $\blamb\in \Pcal_r$, let $\Delta_\blamb(\sigma_+):=\det(\sigma_{\lambda_j+j-i})_{1\leq i,j\leq r}$. Then Giambelli's formula holds:
$$
\wb_{m+\blamb}=\Delta_\blamb(\sigma_+)\wb_m.
$$
}
\eclm
\proof 
One has, by virtue of Proposition \ref{prop:prop1}:
$$
\Delta_\blamb(\sigma_+)\wb_m=\Delta_\blamb(\sigma_+)\wb_m^r\w \wb_{m-r}.
$$
Now $\Delta_\blamb(\sigma_+)\wb_m^r$ can be read in the exterior power $\bw^rM_{\geq m-r+1}$, and then one may invoke \cite[Formula (17)]{SCHSD} or \cite[Corollary 5.8.2]{GatSal4}.\qed

\claim{\bf Boson--Fermion Correspondence.} \label{sec514:bfc} Let $\zeta\in \End_\ZZ(F)$ defined by $\zeta\wb_{m+\blamb}=\wb_{m+1+\blamb}$. It is an automorphism of $F$ and may be thought of as an extension to $F_m$ of the determinant of the map $(\sigma_1)_{|M}$. Let now $B:=\ZZ[e_1,e_2,\ldots]$ be the polynomial ring in the infinitely many indeterminates $(e_1,e_2,\ldots)$ with $\ZZ$-coefficients. Let $E(z):=1+\sum_{j\geq 1}(-1)^1e_jz^j$ and  $H:=(h_j)_{j\in\ZZ}$ be the sequence in $B$ defined by the equality
\be
H(z):=\sum_{j\in \ZZ}h_jz^j:={1\over E(z)}.\label{eq:ezhz}
\ee
From~\eqref{eq:ezhz} it turns out that $h_j=0$ for $j<0$, $h_0=1$ and, for $j>0$, $h_j$ is a $\ZZ$-polynomial in the $e_i$s, homogeneous of degree $j$,  once each $e_i$ is given degree $i$. It is well known that  \cite{MacDonald}
$$
B:=\bigoplus_{\blamb\in\Pcal}\ZZ\Delta_\blamb(H),
$$
where
$$
\Delta_\blamb(H):=\det(h_{\lambda_j-j+i})_{1\leq i,j\leq r},
$$
and $r$ is any positive integer bigger or equal than $\ell(\blamb)$.

This enables to equip $F_0$ with a $B$-module structure by declaring that $h_i\bfu=\sigma_i(\bfu)$ for all $\bfu\in F_0$. In particular 
\be
\sigma_+(z)\wb_{m+\blamb}={1\over E(z)}\wb_{m+\blamb},\label{eq5:bmodsF}
\ee
for all $\blamb\in\Pcal$. 
Since $F_m=\zeta^mF_0$, each of them inherits a structure of free $B$-module generated by $\zeta^m\wb_0$. Consider the polynomial ring $B[\zeta,\zeta^{-1}]$.

\bclm{\bf Definition.} \label{def:bfcorr} {\em 
The Boson--Fermion correspondence is the $B[\zeta,\zeta^{-1}]$-module structure of $F$
$$
B[\zeta,\zeta^{-1}]\otimes F\sra F,
$$
defined by
\be
\Delta_\blamb(H)\wb_m=\wb_{m+\blamb}=\Delta_\blamb(H)\bfb^r_{m+\blamb}\w \wb_{m-r}.\label{eq:bfcor}
\ee
}
\eclm
In particular $\Delta_\blamb(H)\wb_0=\wb_{0+\blamb}$ and 
$
\zeta^m\Delta_\blamb(H)\wb_0=\Delta_\blamb(H)\zeta^m\wb_0=\Delta_\blamb(H)\wb_m.
$
Notice that, on the dual side, i.e. in $F(M^*)$, one has
$\zeta^m\wbet_0=\wbet_{-m}$. In fact 
$$
\delta_{m-n,0}=\delta_{m,n}=\wbet_m(\wb_n)=\wbet_m(\zeta^n\wb_0)=\zeta^n\wbet_m(\wb_0
$$
from which $\zeta^n\wbet_m=\wbet_{m-n}$.

\section{Vertex Operators}\label{sec:sec6}

Via the Boson--Fermion correspondence the Schubert derivations $\ovsig_-(z),\sigma_-(z)$ induce natural maps $B\sra B[z^{-1}]$.

\bclm{\bf Definition.} {\em Define 
$
\ovsig_-(z),\sigma_-(z):B\sra B[z^{-1}]
$
through the equalities
$$
\left(\ovsig_-(z)\Delta_\blamb(H)\right)\wb_0:=\ovsig_-(z)\left(\Delta_\blamb(H)\wb_0\right)=\ovsig_-(z)\wb_{0+\blamb}
$$
and 
$$
(\sigma_-(z)\Delta_\blamb(H))\wb_0:=\sigma_-(z)\left(\Delta_\blamb(H)\wb_0\right)=\sigma_-(z)\wb_{0+\blamb}.
$$
}
\eclm
\bclm{\bf Proposition.}
\begin{eqnarray}
\ovsig_-(z)h_n&=&h_n-h_{n-1}z^{-1}, \label{eq:ovs-hn}\\ \cr \sigma_-(z)h_n&=&\sum_{j\geq 0}{h_{n-j}z^{-j}}.\label{eq:s-hn}
\end{eqnarray}
\eclm

\proof
Recall that $h_j=0$ for $j<0$.
Let us prove \eqref{eq:ovs-hn} first.
\begin{eqnarray*}
(\ovsig_-(z)h_n)\wb_0&=&\ovsig_-(z)(b_n\w\wb_{-1})\cr \cr &=&(b_n-b_{n-1}z^{-1})\w \wb_{-1}=(h_n-h_{n-1}z^{-1})\wb_{0},
\end{eqnarray*}
whence \eqref{eq:ovs-hn}.
The proof of \eqref{eq:s-hn} is analogous:
\begin{eqnarray*}
(\sigma_-(z)h_n)\wb_0&=&\sigma_-(z)(h_n\wb_0)=\sigma_-(z)\sigma_n\wb_0\cr\cr 
&=& \sigma_-(z)(b_n\w \wb_{-1})=\sigma_-(z)b_n\w \wb_{-1}\cr\cr
&=&\sum_{j\geq 0} b_{n-j}\w \wb_{-1}z^{-j}=\sum_{j\geq 0}h_{n-j}z^{-j}\wb_0
\end{eqnarray*}
whence \eqref{eq:s-hn}. \qed

\bclm{\bf Remark.} The sum \eqref{eq:s-hn} is an infinite sum but its multiplication by $\wb_m$ is finite, for all $m\in\ZZ$. For instance $\sum_{j\geq 0}h_{n-j}z^{n-j}\wb_0=\sum_{0\leq j\leq n}h_{n-j}z^{n-j}\wb_0$.
\eclm

Let $\sigma_-(z)H$ denote the sequence $(\sigma_-(z)h_n)_{n\in\ZZ}$ (respectively $\ovsig_-(z)H=(\ovsig_-(z)h_n)_{n\in\ZZ}$).
The following is one of the main results concerning the combinatorics of the subject.
\bclm{\bf Theorem.} {\em Schur determinants commute with taking $\ovsig_{-}(z)$:
\be 
\ovsig_-(z)\Delta_\blamb(H)=\Delta_\blamb(\ovsig_-(z)H)\qquad\mathrm{and}\qquad \sigma_-(z)\Delta_\blamb(H)=\Delta_\blamb(\sigma_-(z)H),
\ee
}
\eclm
\proof
By Proposition \ref{prop:prop1}, the equality
$\wb_{0+\blamb}=\Delta_\blamb(H_r)\wb_0$ holds in the exterior power $\bw^rM_{\geq -r+1}$. We contend that
$$
(\ovsig_-(z)\Delta_\blamb(H_r))\wb_0=\Delta_\blamb(\ovsig_-(z)H)\wb_0,
$$
 and this is true by \cite[Theorem 5.7]{pluckercone} that relies on a general determinantal formula in a polynomial ring due to Laksov and Thorup \cite[Theorem 0.1]{LakTh04}. The same argument holds verbatim for $\sigma_-(z)$. 
 
 \qed
\bclm{\bf Corollary.} {\em $\ovsig_-(z),\sigma_-(z):B\sra B[z^{-1}]$ are rings homomorphisms. 
}
\eclm
\proof In fact $B=\ZZ[h_1,h_2,\ldots]$. Then
\begin{eqnarray*}
\sigma_-(z)(h_{i_1}\cdots h_{i_r})&=&\sigma_-(z)(\sum_{\lambda}a_\blamb\Delta_\blamb(H))\cr\cr &=&\sum_\blamb a_\lambda\Delta_\blamb(\sigma_-(z)H)=\sigma_-(z)h_{i_1}\cdots\sigma_-(z)h_{i_r}.
\end{eqnarray*}
the proof for  $\ovsig_-(z)$ is totally analogous.\qed

In other words, $\sigma_-(z), \ovsig_-(z):B\sra B[[z^{-1}]]$ are Hasse-Schmidt derivations on $B$, in the genuine sense of e.g. \cite[p.~207]{MatsumuraCR}.
\bclm{\bf Lemma.} {\em For all $m\in\ZZ$, $r\in \NN$ and $\blamb\in\Pcal_r$:
\be
\bfb^r_{m+\blamb+(1^r)}\w b_{m-r+1}= \wb^{r+1}_{m+1+\blamb}.\label{eq:12}
\ee
}
\eclm
\proof
The definition of the left--hand side of (\ref{eq:12}) is
$$ 
b_{m+1+\lambda_1}\w\cdots\w b_{m+1-r+1+\lambda_r}\w b_{m+r-1}
$$
$$
=b_{m+1+\lambda_1}\w\cdots\w b_{m+1-(r+1)+2+\lambda_r}\w b_{m+1-(r+1)+1}
$$
that is precisely the definition of its right--hand side.\qed
\bclm{\bf Corollary.}
$$
\bfb^r_{m+\blamb+(1^r)}\w \wb_{m-r+1}=\wb_{m+1+\blamb}.
$$
\eclm
\proof
In fact
\begin{eqnarray*}
\hskip60pt \bfb^r_{m+\blamb+(1^r)}\w \wb_{m-r+1}&=&\bfb^r_{m+\blamb+(1^r)}\w b_{m-r+1}\w\wb_{m-r}\cr\cr
&=&\wb^{r+1}_{m+1+\blamb}\w \wb_{m-r}=\wb_{m+1+\blamb}.\hskip110pt \qed
\end{eqnarray*}
\claim{}\label{sec:rz} Let $R(z)\in \End_{\ZZ[z,z^{-1}]}(F[[z,z^{-1}])$ given by $R(z)\wb_{m+
\blamb}=z^{m+1}\wb_{m+1+\blamb}$. It is clearly invertible: $R(z)^{-1}\wb_{m+
\blamb}=z^{-m}\wb_{m-1+\blamb}$, i.e. $R(z)^{-1}$ is an endomorpism of $F$ homogeneous of 
degree $-1$. On the bosonic side,  define $(R(z)\zeta^m\Delta_\blamb(H))\wb_0=z^{m
+1}\wb_{m+1+\blamb}$, so  that $R(z)\zeta^m\Delta_\blamb(H)=z^{m+1}\zeta^{m+1}
\Delta_\blamb(H)$.
It is a homogeneous operator of degree $1$ (and is the fermionic counterpart of the $R$ operator mentioned in~\cite[Theorem 5.1]{KacRaRoz}).
Notice that for all $j\in\ZZ$, one has
$
R(z)(\sigma_j\wb_{m+\blamb})=\sigma_j(R(z)\wb_{m+\blamb}),
$
as an immediate check shows, whence the commutativity rules
\be
\ovsig_\pm(z)R(z)=R(z)\ovsig_{\pm(z)}\qquad \mathrm{and}\qquad \sigma_{\pm}(z)R(z)=R(z)\sigma_{\pm}(z).\label{eq6:69}
\ee

\bclm{\bf Proposition.}\label{propo:1stvo} {\em Let $\bfb(z):=\sum_{i\in\ZZ}b_iz^{i}\in M[[z^{-1},z]]$. Then for all $\phi\in F$
$$
\bfb(z)\w \phi=R(z)\sigma_+(z)\ovsig_-(z)\phi
$$
}
\eclm
\proof
Each $\phi\in F$ is an integral finite linear combination of homogeneous elements, each summand belonging to $F_m$ for some $m$. Moreover each element of $F_m$ is an integral linear combination of basis elements of the form $\wb_{m+\blamb}$  of $F_m$. Hence it suffices to prove the proposition for $\phi=\wb_{m+\blamb}$.
\begin{center}

\medskip

\begin{tabular}{rcll}
&&$\bfb(z)\w \wb_{m+\blamb}$&\cr\cr
&$=$&$ \displaystyle{\sum_{i\in\ZZ}}b_iz^i\w \bfb^r_{m+\blamb}\w \wb_{m-r}$&\hskip-19pt(definition of $\bfb(z)$  and $\wb_{m+\blamb}$)\cr\cr
&$=$&$\displaystyle{\sum_{i\geq m-r+1}}b_iz^i\w \bfb^r_{m+\blamb}\w \wb_{m-r}$
&\hskip-98pt{\small(the wedge products vanish for $i<m-r+1$)}\cr\cr
&$=$&$z^{m-r+1}\sigma_+(z)b_{m-r+1}\w \bfb^r_{m+\blamb}\w \wb_{m-r}$& \hskip12pt (by definition of $\sigma_+(z)$)\cr\cr
&$=$&$z^{m-r+1}\sigma_+(z)(b_{m-r+1}\w\ovsig_+(z)\bfb^r_{m+\blamb}\w \ovsig_+(z)\wb_{m-r})$&\hskip5pt(integration by parts \eqref{eq:intprt1})\cr\cr
&$=$&$z^{m-r+1}\sigma_+(z)(z^r\ovsig_-(z)\wb^{r}_{m+\blamb+(1^r)}\w b_{m-r+1}\w\ovsig_+(z)\wb_{m-r}$)&\hskip27pt (Definition of $\ovsig_-(z)$)
\end{tabular}

\smallskip
\medskip
\begin{tabular}{rcll}
&$=$&$z^{m+1}\sigma_+(z)(\ovsig_-(z)\wb^{r}_{m+\blamb+(1^r)}\w\wb_{m-r-1})$&\hskip93pt ($z^{m+1}=z^{m-r+1}z^r$)\cr\cr
&$=$&$z^{m+1}\sigma_+(z)(\ovsig_-(z)\wb^{r}_{m+\blamb+(1^r)}\w\ovsig_-(z)\wb_{m-r-1})$&\hskip14pt (the $\ovsig_-(z)$ action (\ref{eq5:sigact}) on $\wb_{m-r+1}$)\cr\cr
&$=$&$z^{m+1}\sigma_+(z)\ovsig_-(z)(\wb^{r}_{m+\blamb+(1^r)}\w\wb_{m-r-1})$&\hskip62pt($\ovsig_-(z)$ is a HS derivation)\cr\cr
&$=$&$z^{m+1}\sigma_+(z)\ovsig_-(z)\wb_{m+1+\blamb}$\cr\cr
 &$=$&$R(z)\sigma_+(z)\ovsig_-(z)\wb_{m+\blamb}$&\hskip40pt(by definition of the map  $R(z)$)\cr\cr
&$=$&$R(z)\cdot\displaystyle{1\over E(z)}\ovsig_-(z)\wb_{m+\blamb}$&{\small(by  the $B$-module structure (\ref{eq5:bmodsF}) of $F$).} \qed
\end{tabular}
\end{center}

\bclm{\bf Corollary.} {\em Let $\Gamma(z):B[\zeta,\zeta^{-1}]\sra B[\zeta,\zeta^{-1}][[z]]$ be given by
$$
(\Gamma(z)\zeta^m\Delta_\blamb(H))\wb_0:=\bfb(z)\w \wb_{m+\blamb}.
$$
Then 
$$
\Gamma(z)=R(z)\sigma_+(z)\ovsig_-(z)=R(z){1\over E(z)}\ovsig_-(z).
$$

}
\eclm
\proof Obvious from the definition. \qed

\bclm{\bf Lemma.} {\em
Let $R(z)^T:F^*[[z,z^{-1}]\sra F^*[[z,z^{-1}]$ be the transpose of the operator\linebreak  $R(z):F[[z,z^{-1}]]\sra F[[z,z^{-1}]]$ as in~{\em \ref{sec:rz}}. Then
\be
R(z)^T\wbet_{m+\blamb}=z^{m}\wbet_{m-1+\blamb},\label{eq:RTr}
\ee
and
\be
(R(z)^T)^{-1}\wbet_{m-1+\blamb}=(R(z)^{-1})^T\wbet_{m-1+\blamb}=z^{-m}\wbet_{m+\blamb}.\label{eq:RTrinv}
\ee
}
\eclm
\proof To prove (\ref{eq:RTr}):
\begin{eqnarray*}
(R(z)^T\wbet_{m+\bmu})(\wb_{m-1+\blamb})&=&\wbet_{m+\bmu}(R(z)\wb_{m-1+\blamb})\cr\cr
&=&\wbet_{m+\bmu}(z^m\wb_{m+\blamb})=z^m\delta_{\bmu,\blamb},
\end{eqnarray*}
and then $R(z)\wbet_{m+\bmu}=z^m\wbet_{m-1+\bmu}$ as stated in~(\ref{eq:RTr}). The proof of (\ref{eq:RTrinv}) is straightforward.\qed

\bclm{\bf Corollary.} \label{cor:prevg*}  {\em If ${\bm\beta}(z)=\sum_{j\in\ZZ} \beta_jz^{-j-1}\in M^*[[z^{-1},z]]$, then
\begin{eqnarray}
{\bm\beta}(z)\w \wbet_{m-1+\bmu}&=&z^{-1}(R(z)^{-1})^T\sigma^T_-(z)\ovsig^T_+(z)\wbet_{m+\bmu}.\label{eq6:Rzt}
\end{eqnarray}
}
\eclm
\proof
Let $\ovsig^*_+(w),\sigma^*_+(w)$ be the Schubert derivations on $\wM^*$ as in Remark \ref{rmk:rmk1}. Let $\beta(w)=w\sum_{j\in\ZZ}\beta_jw^j$. Then applying Proposition \ref{propo:1stvo}
\begin{eqnarray*}
{\bm\beta}(w)\w\wbet_{m-1+\bmu}&=&w\sum_{j\in\ZZ}\beta_jw^{j}\w\wbet_{m-1+\mu}\cr\cr
&=&wR(w)\sigma^*_+(w)\ovsig_-^*(w)\wbet_{m-1+\bmu}=wR(w)^T\sigma^*_+(w)\ovsig_-^*(w)\wbet_{m+\bmu}.
\end{eqnarray*}
Putting $w=z^{-1}$ and observing that ${\ovsig^*_-(w)\beta_j}_{|w=z^{-1}}=\ovsig_+(z)^T\beta_j$ and ${\sigma^*_+(w)\beta_j}_{|w=z^{-1}}=\sigma_-(z)^T\beta_j$, for all $j\in\ZZ$, one finally obtains (\ref{eq6:Rzt}).\qed

\bclm{\bf Proposition.}{\em
\begin{eqnarray}
{\bm\beta}(z)\lrcorner\wb_{m+\blamb}&=&z^{-1}R(z)^{-1}\ovsig_+(z)\sigma_-(z)\wb_{m+\blamb}\cr\cr &=&z^{-1}R(z)^{-1}E(z)\sigma_-(z)\wb_{m+\blamb}.\label{eq:preg*}
\end{eqnarray}
where the last equality follows from (\ref{eq5:bmodsF})}.

\eclm
\proof
For all $\blamb,\bmu\in \Pcal\times\Pcal$ and all $m\in\ZZ$:
$$
\wbet_{m-1+\bmu}({\bm \beta}(z)\lrcorner \wb_{m+\blamb})=({\bm\beta}(z)\w \wbet_{m-1+\bmu})(\wb_{m+\blamb}).
$$
By Corollary \ref{cor:prevg*}, then
\begin{eqnarray*}
\wbet_{m-1+\bmu}({\bm \beta}(z)\lrcorner \wb_{m+\blamb})&=&z^{-1}(R(z)^{-1})^T\sigma^T_-(z)\ovsig^T_+(z)\wbet_{m+\bmu}\big(\wb_{m+\blamb}\big)\cr\cr
&=&\wbet_{m-1+\bmu}(z^{-1}\ovsig_+(z)\sigma_-(z)R(z)^{-1}\wb_{m+\blamb}),
\end{eqnarray*}
whence \eqref{eq:preg*}, because of \eqref{eq6:69}.\qed

\bclm{\bf Corollary.} {\em Let $\Gamma^*(z):=B[\zeta,\zeta^{-1}]\sra B[\zeta,\zeta^{-1}][[z]]$ be given by
$$
(\Gamma^*(z)\zeta^m\Delta_\blamb(H))\wb_0:={\bm\beta}(z)\lrcorner \wb_{m+\blamb}.
$$
Then 
$$
\Gamma^*(z)=z^{-1}R(z)^{-1}\ovsig_+(z)\sigma_-(z)=z^{-1}R(z)^{-1}E(z)\sigma_-(z).
$$

}
\eclm
\proof A straightforward consequence of the definition.\qed
\section{Another expression for $\Gamma^*(z)$}\label{sec:sec7}
\claim{} Notation as in Section \ref{sec:nota}.
Expanding the determinant of the Schur polynomial $\Delta_\blamb(H)$ along the first row, according to Laplace's rule, one easily check that
\be 
\Delta_\blamb(H)=h_{\lambda_1}\Delta_{\blamb^{(1)}(H)}-h_{\lambda_2-1}\Delta_{\blamb^{{(2)}}+(1)}(H)+\cdots+(-1)^{r-1}h_{\lambda_r-r+1}\Delta_{\blamb^{(r)}+(1^{r-1})}(H).
\ee
Unable to find a better compact notation for it, we denote by $\Delta_\blamb(z^{-\blamb},H)$ the determinant obtained by $\Delta_\blamb(H)$ under the substitution $h_{\lambda_j-i+1}\sra z^{-\lambda_j+i-1}$. In other words:
\be
\Delta_\blamb(z^{-\blamb},H):={\Delta_{\blamb^{(1)}}(H)\over z^{-\lambda_1}}-{\Delta_{\blamb^{(2)}}(H)\over z^{\lambda_2-1}}+\cdots+(-1)^{r-1}{\Delta_{\blamb^{(r)}+(1^{r-1})}(H)\over z^{\lambda_r-r+1}}
\ee
or, more explicitly
\be 
\Delta_\blamb(z^{-\blamb},H):=\left|\begin{matrix}\displaystyle{1\over z^{\lambda_1}}&\displaystyle{1\over z^{\lambda_2}+1}&\cdots& \displaystyle{1\over z^{\lambda_r+r-1}}\cr\cr
h_{\lambda_1+1}&h_{\lambda_2}&\cdots&h_{\lambda_r+r-2}\cr
\vdots&\vdots&\ddots&\vdots\cr
h_{\lambda_1+r-1}&h_{\lambda_2+r-2}&\cdots&h_{\lambda_r}\end{matrix}\right|.
\ee
Then, we have:
\bclm{\bf Proposition.}
\be
{\bm\beta}(z^{-1})\lrcorner \wb_{m+\blamb}=z^{-m-1}\left(\Delta_\blamb(z,H)+(-1)^{r-1}\sum_{j\geq 0}(-1)^j\Delta_{\blamb+(1^{r+j})}(H)z^{j+r}\right)\wb_{m-1}.\label{eq:explvo}
\ee
\eclm
\proof
First of all we apply directly Leibniz rule enjoyed by the derivation ${\bm\beta}(z)\lrcorner$:
\begin{eqnarray}
{\bm\beta}(z^{-1})\lrcorner \wb_{m+\blamb}&=&{\bm\beta}(z^{-1})\lrcorner (\bfb^r_{m+\blamb}\w\wb_{m-r})\cr\cr
&=&{\bm\beta}(z^{-1})\lrcorner \bfb^r_{m+\blamb}\w \wb_{m-r}+(-1)^{r-1}\bfb^r_{m+\blamb}\w\left({\bm\beta}(z^{-1})\lrcorner\wb_{m-r}\right).\cr&&\label{eq:lastlin}
\end{eqnarray}
We compute separately the two summands occurring in the r.h.s. of~\eqref{eq:lastlin}. Let us begin with the second one:
\begin{eqnarray}
\bfb^r_{m+\blamb}\w\left({\bm\beta}(z^{-1})\lrcorner\wb_{m-r}\right)&=&\bfb^r_{m+\blamb}\w z^{-m+r-1}\sigma_+^T(z)\beta_{m-r}\lrcorner \wb_{m-r}\cr\cr
&=&z^{-m-1+r}\bfb^r_{m+\blamb}\w \ovsig_+(z)(\beta_{m-r}\lrcorner \wb_{m-r})\cr\cr
&=&z^{m-1+r}\bfb^r_{m+\blamb}\w \ovsig_+(z)(\beta_{m-r}\lrcorner(\sigma_+(z)b_{m-r}\w \wb_{-m-1+r}))\cr\cr
&=&z^{-m-1+r}\bfb^r_{m+\blamb}\w \ovsig_+(z)\wb_{m-r-1}\cr\cr
&=&z^{-m-1+r}\bfb^r_{m+\blamb}\w\sum_{j\geq 0}(-1)^j\wb_{m-r-1+(1^j)}z^j\cr\cr
&=&z^{-m-1+r}\sum_{j\geq 0}(-1)^j\bfb^r_{m-1+\blamb+(1^r)}\w \wb_{m-1-r+(1^j)}z^j\cr\cr
&=&z^{-m-1}\sum_{j\geq 0}(-1)^j\wb_{m-1+\blamb+(1^{r+j})}.\label{eq:1stsum}
\end{eqnarray}

To compute the first summand, instead,  it is sufficient to apply the definition of contraction: each $b_i$ occurring in the expression $\bfb^r_{m+\blamb}$ is replaced, with the appropriate sign, by $z^{-i-1}$.
The straightforward equality, 
$$
\bfb^r_{m+\blamb}=(-1)^{j}b_{m-j+\lambda_{j+1}}\w \bfb^r_{m+\blamb^{(j)}+(1^{j-1})},\qquad \qquad 1\leq j\leq r,
$$
holding by the very meaning \eqref{eq:notwbr} of the notation $\bfb^r_{m+\blamb}$, easily implies that
\begin{eqnarray}
{\bm\beta}(z^{-1})\lrcorner \bfb^r_{m+\blamb}&=&{1\over z^{m+1}}{\bfb^r_{m+\blamb^{(1)}}\over z^{\lambda_1}}-{\bfb^r_{m-1+(\blamb^{(2)}+(1))}\over z^{\lambda_2-1}}+\cdots+(-1)^{r-1}{\bfb^r_{m-r+1+(\blamb^{(r)}+(1^{r-1}))}\over z^{\lambda_r-r+1}}\cr\cr
&=&{1\over z^{m+1}}\Delta_\blamb(z^{-\blamb},H)\wb_{m-1},\label{eq:secondf}
\end{eqnarray}
where in the last equality we used the Boson--Fermion correspondence \eqref{eq:bfcor}.
Formula \eqref{eq:explvo} then follows by plugging \eqref{eq:1stsum} and \eqref{eq:secondf} into \eqref{eq:lastlin}, and using again the Boson--Fermion correspondence.\qed

\section{\bf Commutation Rules}\label{sec:sec8}

\bclm{\bf Proposition.}\label{prop:precom}
{\em
\begin{enumerate}
\item[i)] The following commutation rules hold in $\End_\ZZ(\bw M)$ for all $i,j\in\ZZ$:
\be 
\sigma_i\sigma_j=\sigma_j\sigma_i,\qquad \ovsig_i\sigma_j=\sigma_j\ovsig_i,\qquad \ovsig_i\ovsig_j=\ovsig_j\ovsig_i.\label{eq1:comrul}
\ee
\item[ii)] All the operators $\sigma_{\pm}(z), \ovsig_{\pm}(z):\bw M\sra \bw M[[z]]$ are mutually commuting.
\end{enumerate}
}
\eclm
\proof
i) In fact for all $k\in\ZZ$, one has $\sigma_i\sigma_jb_k=b_{k+i+j}=\sigma_j\sigma_ib_k$. Assume commutation holds for $\bw^{r-1}M$. Each $u\in \bw^rM$  is a sum  of typical elements of the form $b\w v$ for $b\in M$ and $v\in \bw^{r-1}M$. Then
\begin{eqnarray*}
\sigma_i\sigma_j(b\w v)&=&\sigma_i\sum_{k=0}^j\sigma_kb\w \sigma_{j-k}v\\
&=&\sum_{l=0}^i\sum_{k=0}^j\sigma_l\sigma_kb\w \sigma_{i-l}\sigma_{j-k}v\\
&=&\sum_{l=0}^i\sum_{k=0}^j\sigma_k\sigma_lb\w \sigma_{j-k}\sigma_{i-l}v=\sigma_j\sigma_i(b\w v),
\end{eqnarray*}
where the third equality follows by induction, which proves i). Item ii) is  a straightforward consequence of  i).
\qed

The commutation rules \eqref{eq1:comrul} do not hold when $\sigma_{\pm}(z)$ and $\ovsig_{\pm}(z)$ are extended to the fermionic space $F$. For example,
$$
\sigma_{-1}\sigma_2\wb_0=\sigma_{-1}b_2\w\wb_{-1}=b_1\w\wb_{-1}\neq 0=\sigma_2\cdot 0=\sigma_2\sigma_{-1}\wb_{-1}.
$$
\claim{\bf Notation.} Let $\Rcal(z,w)$ be a rational expression having poles in $z=0$, $w=0$ and $z-w=0$. Following \cite[p.~18]{Kacbeg} we denote by $i_{w,z}R(z,w)$ the expansion of $\Rcal(z,w)$ as a formal power series in $(z/w)$ and by $i_{z,w}$ the expansion of the same expression as a formal power series of $w/z$.
\bclm{\bf Lemma.}\label{lem:precom} {\em For all $m\in\ZZ$, 
$$
\sigma_-(w)\sigma_+(z)\wb_m=i_{w,z}{w\over w-z}\sigma_+(z)\sigma_-(w)\wb_m.
$$
}
\eclm
\proof
\begin{center}
\begin{tabular}{rclr}
$\sigma_-(w)\sigma_+(z)\wb_m$&$=$&$\sigma_-(w)(\sigma_+(z)b_m\w \wb_{m-1})$\hskip 40pt(decomposition of $\wb_m$)&\cr\cr
&$=$&$\sigma_-(w)\sigma_+(z)b_m\w \wb_{m-1}$\hskip 49pt(distributing $\sigma_-(w)$ with&\cr
&&\hskip 180pt respect to $\wedge$)&\cr\cr
&$=$&$\sigma_-(w)\sum_{j\geq 0}b_{m+j}z^j\w \wb_{m-1}$ \hskip 24pt (definition of $\sigma_+(z)$)&\cr\cr
\end{tabular}
\end{center}
\begin{center}
\begin{tabular}{rclr}

&$=$&$\displaystyle{\sum_{j\geq 0}}\displaystyle{\sum_{i=0}^{j}}\displaystyle{b_{m+j}z^j\over w^i} \w \wb_{m-1}$\hskip20pt &\cr\cr
&$=$&$\left(1+\displaystyle{z\over w}+\displaystyle{z^2\over w^2}+{z^3\over w^3}+\cdots\right)\sigma_+(z)\sigma_-(w)b_m \w \wb_{m-1}$\cr\cr
&$=$&$i_{w,z}\displaystyle{w\over w-z}\sigma_+(z)\sigma_-(w)\cdot \wb_m.$\hskip 220pt \qed

\end{tabular}
\end{center}

\bclm{\bf Proposition.} {\em Let $\sigma_{\pm}(z)$ and $\ovsig_{\pm}(z^{\pm 1})$ considered as maps from $F\sra F[[z^{\pm1}]]$. Then the following commutation rules holds
\be 
\sigma_-(w)\sigma_+(z)=i_{w,z}{w\over w-z}\cdot\sigma_+(z)\sigma_-(w)\label{eq:comm1}
\ee

\be 
\ovsig_-(z)\ovsig_+(w)=i_{z,w}{z\over z-w}\ovsig_+(w)\ovsig_-(z)\label{eq:comm2}
\ee
}
\eclm
\proof Let us prove \eqref{eq:comm1}.
For all $\blamb\in\Pcal$ we have
$$
\sigma_-(w)\sigma_+(z)\wb_{m+\blamb}=\sigma_-(z)\sigma_+(w)(\bfb^r_{m+\blamb}\w \wb_{m-r})
$$
provided that $r$ is bigger or equal than the length of the partition $\blamb$. Now we use the definition of the extension of $\sigma_{\pm}(z)$ to $F$. This gives:
\begin{eqnarray*}
&&\sigma_-(w)(\sigma_+(z)\bfb^r_{m+\blamb}\w \sigma_+(z)b_{m-r}\w \wb_{m-r-1})=\cr\cr
&=&\sigma_-(w)\sigma_+(z)\bfb^r_{m+\blamb}\w \sigma_-(w)\sigma_+(z)b_{m-r}\w \wb_{m-r-1}
\end{eqnarray*}
Proposition \ref{prop:precom} ensures that $\sigma_-(w)$ and $\sigma_+(z)$ commute on $\wM$ so obtaining
$$
\sigma_+(z)\sigma_-(w)\bfb^r_{m+\blamb}\w \sigma_-(w)\sigma_+(z)\wb_{m-r}
$$
and now we apply Lemma~\ref{lem:precom} to finally obtain:
\begin{eqnarray*}
\sigma_+(z)\sigma_-(w)\bfb^r_{m+\blamb}\w \sigma_-(w)\sigma_+(z)\wb_{m-r}&=&\sigma_+(z)\sigma_-(w)\left(\bfb^r_{m+\blamb}\w {w\over w-z}\wb_{m-r}\right)\cr
&=&i_{w,z}{w\over w-z}\sigma_+(z)\sigma_-(w)\wb_{m+\blamb}
\end{eqnarray*}
which proves the commutation formula~\eqref{eq:comm1}. To prove \eqref{eq:comm2} we take the inverse of either side of \eqref{eq:comm1} obtaining
$$
\ovsig_+(z)\ovsig_-(w)=\left(1-{z\over w}\right)\ovsig_-(w)\ovsig_+(z),
$$
from which
$$
\ovsig_-(w)\ovsig_+(z)=i_{w,z}{w\over w-z}\ovsig_+(z)\ovsig_-(w).
$$
Changing the role of the indeterminates $z$ and $w$ one obtains precisely \eqref{eq:comm2}.\qed

\section{The DJKM Bosonic Vertex Representation of $\glin$}\label{sec:sec9}

\claim{} Let $\Bcal_{ij}:=b_i\otimes\beta_j\in M\otimes M^*\cong End(M)$ and $gl_\infty(M):=\bigoplus_{i,j\in\ZZ}\ZZ \Bcal_{ij}$. It is a Lie algebra with respect to the obvious Lie bracket $[A,B]=AB-BA$ ($A,B\in \glin$). Let
$$
\delta:\glin\sra End_\ZZ(\wM)
$$
be the representation \eqref{eq0:2} of $\glin$ as a sub-algebra of derivations of $\wM$.
%\begin{eqnarray*}
%\delta(A)\bfu&=&A\bfu, \qquad \forall \bfu\in M,\cr\cr
%\delta(A)(\bfu\w \bfv)&=&\delta(A)u\w v+u\w \delta(A)\bfv,\qquad \forall \bfu,\bfv\in\wM.
%\end{eqnarray*}
\bclm{\bf Proposition.} {\em Let $b\otimes\beta\in M\otimes M^*$ and $\bfu\in\wM$. Then 
$$
\delta(b\otimes\beta)(\bfu)=b\w (\beta\lrcorner \bfu).
$$
}
\eclm
\proof As $\bfu\in \wM$ is a finite sum of homogeneous elements,  we may assume without loss of generality that $u\in \bw^rM$. Then we argue on induction on $r\geq 1$. If $\bfu\in M$, $\delta(b\otimes \beta)(\bfu)=\beta(\bfu)b=b(\beta\lrcorner \bfu)$, and the claim holds for $r=1$. Assume now the property true for all $\bfu\in\bw^iM$ and  $1\leq i\leq r-1$. Each $\bfv\in \bw^rM$ is a finite sum of monomials of the form $u\w w$, with $\bfu\in M$ and $\bfw\in \bw^{r-1}M$. We may then assume $v=u\w w$ and, in this case,
\begin{eqnarray*}
\delta(b\otimes\beta)(\bfu\w \bfw)&=&\beta(\bfu)b\w \bfw+\bfu\w b\w (\beta\lrcorner \bfw)\cr\cr
&=&b\w (\beta(\bfu)\bfw-\bfu\w \beta \lrcorner \bfw)\cr\cr &=&b\w (\beta\lrcorner (\bfu\w \bfw)).
\end{eqnarray*}
\qed

We extend the derivation $\delta$ of $\wM$ to $F_m$ as follows. Each $A\in \glin$ is a finite linear combination $\sum_{ij}a_{ij}b_i\otimes \beta_j$. Let $k$ be the minimum among all $j$ such that $a_{ij}\neq 0$ and let  $r\geq 0$ such that $m-r<k$. Thus one defines $\delta_m:\glin\sra End_\ZZ(F_m)$ via:
$$
\delta_m(A)\wb_{m+\blamb}=\delta_m(A)\wb_{m+\blamb}^r\w \wb_{m-r}.
$$
An easy check shows that the definition does not depend on the choice of  $r\geq 0$ such that $m-r<k$.
Let 
\[
\delta_m(z,w)=\sum_{i,j\in\ZZ}\delta_m(\Bcal_{ij})z^iw^{-j}:F_m\sra F_m[[z,w^{-1}]]
\]

\bclm{\bf Theorem} \cite[Date--Jimbo--Kashiwara--Miwa]{DJKM01}
\begin{eqnarray}
\delta_m(z,w)=\sum_{i,j\in\ZZ}\delta_m(\Bcal_{ij})z^iw^{-j}&=&{z^{m}\over w^{m}}\sigma_+(
z)\ovsig_-(z)\ovsig_+(w)\sigma_-(w)\cr 
&=&{z^{m}\over w^{m}}i_{z,w}{z\over z-w}\sigma_+(z)\ovsig_+(w)\ovsig_-(z)\sigma_-(w)\cr\cr
&=&{z^{m}\over w^{m}}i_{z,w}{z\over z-w}{E(w)\over E(z)}\ovsig_-(z)\sigma_-(w).
\label{eq:fermvert}
\end{eqnarray}

\eclm
\proof
We have
\begin{eqnarray*}
\sum\delta_m(\Bcal_{ij})\wb_{m+\blamb}z^iw^{-j}&=&\sum \delta_m(b_i\otimes \beta_j)\wb_{m+\blamb}z^iw^{-j}
\cr\cr
&=&\bfb(z)\w (w{\bm \beta}(w)\lrcorner \wb_{m+\blamb}))
\cr\cr
&=&\bfb(z)\w (w^{-m}\ovsig_+(w)\sigma_-(w)\wb_{m-1+\blamb})\cr\cr
&=&w^{-m}\bfb(z)\w\ovsig_+(w)\sigma_-(w)\wb_{m-1+\blamb}.
\end{eqnarray*}
Now $\ovsig_+(w)\sigma_-(w)\wb_{m-1+\blamb}$ is a $\ZZ[[w,w^{-1}]$-linear combination of  elements of $F_{m-1}$ and $b(z)\w$ is $\ZZ[[w,w^{-1}]$ linear. Proposition \ref{propo:1stvo} applied to $F_{m-1}$ gives

\begin{eqnarray*}
&&z^{m}\sigma_+(z)\ovsig_-(z^{-1})w^{-m}\ovsig_+(w)\sigma_-(w)\wb_{m+\blamb}
\cr\cr
&=&{z^{m}\over w^{m}}\sigma_+(z)\ovsig_-(z)\ovsig_+(w)\sigma_-(w)\wb_{m+\blamb}\cr\cr
&=&{z^{m}\over w^{m}}i_{z,w}{z\over z-w}\sigma_+(z)\ovsig_+(w)\ovsig_-(z^{-1})\ovsig_-(w)\wb_{m+\blamb},
\end{eqnarray*}
and using the definition of the $B$-module structure of $F_m$ one gets precisely (\ref{eq:fermvert}).\qed

Let
\[
(\Bcal^{(m)}(z,w)\Delta_\blamb(H))\wb_m=\delta_m(z,w)\wb_{m+\blamb}.
\]
\bclm{\bf Corollary (The DJKM bosonic Vertex Representation of $\glin$).}
{\em Using the Boson--Fermion correspondence \ref{sec514:bfc}:
$$
\Bcal^{(m)}(z,w)={z^m\over w^m}i_{z,w}{z\over z-w}\Gamma(z,w),
$$
where the {\em vertex operator} $\Gamma(z,w)$ is given by
\be
\Gamma(z,w)={E(w)\over E(z)}\ovsig_-(z)\sigma_-(w).\label{eq:vop}
\ee
}
\eclm
\proof
It follows straightforwardly from the definition and  expression of (\ref{eq:fermvert}) for $\delta_m$.\qed

\claim{} Let $\Acal_\infty(\ZZ)$ be the Lie algebra of matrices $(a_{ij})_{i,j\in\ZZ}$ having only finitely many non-zero diagonals, i.e. $a_{ij}=0$ if $|i-j|>>0$. In this case the representation $\delta_m$ nust be replaced by a modified representation $\widehat{\delta}_m$ \cite[p.~40]{KacRaRoz}:
$$
\left\{\begin{matrix}
\widehat{\delta}_m(\Bcal_{ij})&=&\delta_m(\Bcal_{ij})&\mathrm{if}& i\neq j&\mathrm{or}&i=j>0,\cr\cr
\widehat{\delta}_m(\Bcal_{ii})-{\mathbf 1}_{F_m}&=&\delta_m(\Bcal_{ii})&\mathrm{if}&& &\,\,\,\,\,i=j\leq 0.
\end{matrix}\right. 
$$
Then, to obtain the generating function of the representation of $\Bcal_{ij}$ via $\widehat{\rr}_m$ it suffices to subtract from formula \eqref{eq:fermvert} the series 
$$
\sum_{j\leq 0}z^jw^{-j}\wb_{m+\blamb}=\sum_{j\geq 0} \left({w\over z}\right)^j\wb_{m+\blamb}=i_{z,w}{z\over z-w}\wb_{m+\blamb},
$$
so obtaining
$$
\sum_{i,j\in\ZZ}\widehat{\delta}_m(\Bcal_{ij})z^iw^{-j}=i_{z,w}\,{z\over z-w}\cdot\left({z^m\over w^m}\Gamma(z,w)-1\right),
$$
where $\Gamma(z,w)$ is like in (\ref{eq:vop}).

\claim{\bf Remark.} Recall that $\Gamma(z,w)=\displaystyle{E(w)\over E(z)}\ovsig_-(z)\sigma_-(w)$ is a well defined operator $B\sra B[[z,w^{-1}]$, which is defined over the  integers. In  $B_\QQ:=B\otimes_\ZZ \QQ$ one can define variables $(x_1,x_2,\ldots)$ through the equalities \cite{KacRaRoz}:
\be
\exp\left(\sum_{i\geq 1}x_iz^i\right)={1\over E(z)}\qquad \mathrm{and}\qquad \exp\left(-\sum_{i\geq 1}x_iw^i\right)=E(w).\label{eq:exp1}
\ee
Moreover in \cite[Theorem 5.7]{pluckercone} it is shown that, over the rationals,
\be
\ovsig_-(z)=\exp\left(-\sum_{i\geq 1}{1\over iz^i}{\partial\over \partial x_i}\right)\qquad \mathrm{and}\qquad \sigma_-(w)=\exp\left(\sum_{i\geq 1}{1\over iw^i}{\partial\over \partial x_i}\right),\label{eq:exp2}
\ee
so that after substituting (\ref{eq:exp1}) and (\ref{eq:exp2}) into expression (\ref{eq:vop}) returns the traditional form (\ref{eq0:clasvo}).
%$$
%\Gamma(z,w)\otimes_\ZZ{1_\QQ}=\exp\left(\sum_{i\geq 1}x_i(z^i-w^i)\right)\exp\left(-\sum_{i\geq 1}{1\over i}\left({1\over z^i}-{1\over w^i}\right){\partial\over \partial x_i}\right),
%$$
%(see e.g. \cite[Formula 5.33]{KacRaRoz})
\bibliographystyle{amsplain}
\providecommand{\bysame}{\leavevmode\hbox to3em{\hrulefill}\thinspace}
\providecommand{\MR}{\relax\ifhmode\unskip\space\fi MR }
\providecommand{\MRhref}[2]{%
  \href{http://www.ams.org/mathscinet-getitem?mr=#1}{#2}
}
\providecommand{\href}[2]{#2}

%\bibliography{../files_comuni_tex/mia_biblio}

\parbox[t]{3in}{{\rm Letterio~Gatto}\\
{\tt \href{mailto:letterio.gatto@polito.it}{letterio.gatto@polito.it}}\\
{\it Dipartimento~di~Scienze~Matematiche}\\
{\it Politecnico di Torino}\\
{\it ITALY}} \hspace{1.5cm}
\parbox[t]{2.5in}{{\rm Parham~Salehyan}\\
{\tt \href{mailto:p.salehyan@unesp.br}{p.salehyan@unesp.br}}\\
{\it Ibilce UNESP}\\
{\it Campus de S\~ao Jos\'e do Rio Preto, SP}\\
{\it BRAZIL}}

\end{document}